\documentclass[12pt,twoside]{article}
\usepackage{a4wide}
\usepackage{amsmath}
\usepackage{theorem}
\usepackage{amssymb}
\usepackage[dvips]{graphicx}
\usepackage[dvips]{color}

\theoremstyle{changebreak}

\newtheorem{thm}{Theorem}[section]
\newtheorem{defn}[thm]{Definition}
\newtheorem{lem}[thm]{Lemma}
\newtheorem{cor}[thm]{Corollary}

\newtheorem{rmk}[thm]{Remark}

\def\proof{\par\smallskip
             \noindent {\sc Proof. }}
\def\proofof #1 {\par\medskip\noindent {\sc Proof of #1. }}

\def\sketch{\par\medskip\noindent{\sc Sketch of proof. }}
\def\sketchof #1 {\par\medskip\noindent {\sc Sketch of proof of #1. }}

\def\qed{\rule{0pt}{0pt}\hfill $\blacksquare$ \par\medskip}

\def\remark{\par\medskip \noindent {\sc Remark. }}

\renewcommand{\mod}{\operatorname{mod}\,}

\newcommand{\Lbac}{\backslash}

\newcommand{\N}{\mathbb{N}}

\newcommand{\R}{\mathbb{R}}

\newcommand{\C}{\ensuremath{\mathbb{C}}}
\newcommand{\Ch}{\hat{\mathbb{C}}}

\newcommand{\D}{D}

\newcommand{\OO}{\mathcal{O}}
\newcommand{\oo}{\emph{o}}

\newcommand{\dist}{\operatorname{dist}}
\newcommand{\diam}{\operatorname{diam}}
\newcommand{\meas}{\operatorname{meas}}
\newcommand{\sing}{\operatorname{sing}}

\newcommand{\ep}{\epsilon}

\newcommand{\G}{G}

\newcommand{\F}{{\cal{F}}}
\renewcommand{\S}{{\cal{S}}}

\begin{document}
\sloppy
\title{Recurrence of entire transcendental functions with simple post singular sets}
\author{Jan-Martin Hemke, CAU-Kiel}
\date{\today}
\maketitle  \abstract{\quad\\We study how the orbits of the
singularities of the inverse of a meromorphic function prescribe the
dynamics on its Julia set, at least up to a set of
(Lebesgue) measure zero. We concentrate on a family of entire transcendental functions with only finitely many singularities of the inverse, counting multiplicity, all of which either escape exponentially fast or are pre-periodic. For these functions we are able to decide, whether the function is recurrent of not. In the case that the Julia set is not the entire plane we also obtain estimates for the measure of the Julia set.}
\section{Introduction}
One of the main ideas in complex dynamics is to divide the plane into the Fatou set of points, where the iterates
behave stable, i.e. where they form a normal family, and its complement, the Julia set. By definition the dynamics in
the Fatou set is the easier and understood very well. We are interested in the dynamics of meromorphic functions
on their Julia set. In \cite{bock} H.
Bock proved the following
\begin{thm}[Bock]\label{bockclass} For any non-constant meromorphic function, which
is defined on the whole complex plane, one of the two following
cases holds:
\begin{itemize}
\item[(i)] The Julia set of $f$ is the entire plane and
for all $A\subset\C$ of positive measure, all $m\in\N$ and almost all $z\in\C$ there are infinitely many $n\in\N$ with $f^{mn}(z)\in A$;
\item[(ii)]  almost every forward-orbit in the Julia set accumulates only in
the post-singular set. 
\end{itemize}
\end{thm}
Here the post-singular set denotes the closure
of the union of the forward-orbits of all singularities of the
inverse function, which are the critical and asymptotic values. This result is a generalization of similar results for rational functions, obtained by M. Lyubich \cite{lyubich2} and C. McMullen \cite{mcmullen2}.\\
We introduce some important terms from ergodic theory, which are related to the classification above. A meromorphic function is called \emph{ergodic} (with respect to the Lebesgue-measure), if any invariant set has full
measure or measure zero. It is called \emph{recurrent}, if for every set $A\subset \C$ and almost every point $z\in A$ the
cardinality of the set $A\cap O^+(z)$ is infinite. It is easy to see that (i) implies recurrence and ergodicity. Case  (ii) does not rule out either one of these two in general. If however $P(f)\not =\Ch$, it implies non-recurrence.\\
It is natural to ask for a given function, which case holds. Since a non-empty Fatou set
always implies (ii), one can restrict to the
cases, in which the Julia set consists of the whole complex plane. If the Julia set is not the entire plane, and thus (ii) holds, it would still be interesting to know whether the Julia set has positive measure, since otherwise the statement (ii) would be trivial. \\
In the paper mentioned, H. Bock
gives sufficient conditions for (i): If $f$ is entire, the set of
singularities of the inverse function is finite and all of these are
pre-periodic but not periodic, then (i) is satisfied.
Thus the function $f(z)=2\pi i \exp(z)$ is an example for this first case, in which the post-singular set consists of the only asymptotic value zero and its image $2\pi i$. Other conditions concerning this case are given by L. Keen and J. Kotus \cite{keenkotus}. Conversely it was already shown in 1984, independently by M. Rees \cite{rees} and M. Lyubich \cite{lyubich}, that the function $f(z)=\exp(z)$ is an example for (ii). Here the post-singular set consists of the the closure of the forward-orbit of the only asymptotic value zero, which tends to infinity on the real axis. This result was generalized in \cite{hemke} to functions $f_{\lambda}(z)=\lambda \exp(z)$, if $f_{\lambda}^n(0)$ tends to infinity sufficiently fast. M. Urbanski and A. Zdunik \cite{urbanski:zdunik} even showed that the Hausdorff-dimension of the remaining set is smaller than 2.\\ The difference between the dynamics of $\exp(z)$ and $2\pi i\exp(z)$ is caused by the different behavior of the asymptotic value zero under iteration.
One might hope for a classification of the two cases depending on the behavior of the singularities of the inverse.
As a first approach, we restrict to functions with few
singularities of the inverse, which have simple orbits. As we said, we can neglect all orbits that would imply the existence of a component of the Fatou set, such as periodic critical points or infinite orbits that converge in $\C$. The easiest orbits that remain are pre-periodic or escaping ones. If one considers meromorphic functions with poles, another interesting case is that of singularities, which are mapped eventually onto a pole. This case has been studied by B. Skorulski for the tangent family in \cite{bartok} and for a larger class of functions in his recent thesis.\\
We are interested in conditions ensuring case (ii). In the third chapter we prove the rather technical theorem \ref{asatz}, which provides a set of sufficient conditions for this case. The proof consists of applying the method developed by M. Rees in order to construct a  positive measure set of points, whose iterates show a ``spiral'' type of behavior: They are eventually mapped close to some asymptotic value, then follow its orbit for a certain number of iterates, coming close to infinity, until they are mapped again, and even closer than before, to some asymptotic value etc. These orbits are not dense in $\Ch$, which implies (ii). Since these orbits accumulate at infinity, the only type of components that could possibly intersect this set, are wandering and Baker domains. At least for the various families (e.g. critically finite entire functions), in which those do not occur (see remark \ref{remk}), we also obtain that the Julia set has positive measure. \\ In chapter 4 we consider functions of the type $f(z)=\int_{0}^z
P(t)\exp(Q(t))dt+c$ with polynomials $P$ and $Q$ and $c\in\C$, such that $Q$ is not constant and $P$ not zero. These functions have at most $\deg(Q)$ asymptotic values and $\deg(P)$ critical points.
In the extremal case that all singularities of the inverse are pre-periodic but not periodic, the theorem of H. Bock implies (i). We consider the other extremal case, in which the singularities of the inverse tend to infinity. It turns out that we may neglect the critical values, but have to specify the speed of escape of the asymptotic values. We say that a point $z$ \emph{escapes exponentially}, if $|f^n(z)|\ge
\exp(|f^{n-1}(z)|^{\delta}) $ for some $\delta>0$ and almost all
$n\in\N$. Then theorem \ref{asatz} yields 
\begin{thm}\label{gsatz}
Let $P$ and Q be polynomials with $P$ not zero and $Q$ not constant, $c\in \C$ and \begin{eqnarray*}f(z):=\int_{0}^z
P(t)\exp(Q(t))dt+c.\end{eqnarray*} Suppose that all asymptotic values escape exponentially. Then the Julia set has positive measure and
$\omega(z)\subset P(f)$ for almost all $z\in J(f)$. If $\deg(Q)\ge 3$, then $\meas(F(f))<\infty$.
\end{thm}
Here $\omega(z)$ denotes the $\omega$-limit set that consists of all accumulation points of the sequence $(f^n(z))$.
Conversely one may ask, whether almost every orbit in the Julia set accumulates at every singularity $s$ of $f^{-1}$. It is easy to find examples, for which this is not the case,
if $s$ is a critical value. In order not to accumulate at an asymptotic value $s$, an orbit has to stay out of an entire
sector. In other contexts, sets with this property turned out to have measure zero. Thus one may expect that indeed for almost every point $z\in J(f)$  every asymptotic values $s$ is contained in $\omega(z)$. If however the set of points in the Julia set, whose orbit are bounded, had positive measure,
there would be no reason why these orbits should accumulate at a given asymptotic value. It is not known, whether this can
actually occur, and related to the question, whether the Julia set of a polynomial may have positive measure, which is a well known open question. A positive answer to this question would suggest a negative answer to our initial question also for asymptotic values. Under additional
assumptions on the critical values however the answer is positive. More precisely we get:
\begin{thm}\label{andereinclusion}
Let $f$ be as above and again suppose all its asymptotic values escape exponentially. Suppose that every critical point either also escapes exponentially, is pre-periodic or is contained in an attractive
Fatou-component. Then $\omega(z)=\overline{O^+(A)}$ for almost every point $z\in J(f)$, where $A$ denotes the set of asymptotic values.
\end{thm}
We define the \emph{multiplicity} of an asymptotic value $s$ as the supremum of the set of all natural numbers $n$ with the following property: There exists an $\ep_0>0$ such that for all $\ep<\ep_0$ the set $f^{-1}(B(s,\ep))$ contains at least $n$ unbounded components. Then the functions from above have exactly $\deg(Q)$ asymptotic values and $\deg(P)$ critical points, counting multiplicity, and may even be characterized as those entire transcendental functions with this property. This was proved by G. Elfving in \cite{elfving}. He generalized a method introduced by R. Nevanlinna from \cite{nevanlinna1}, who showed the same in the case $\deg(P)=0$. This method is summarized in 
\cite{nevanlinna}. 
\begin{thm}[Elfving] Let $f$ be entire transcendental, with only finitely many singularities of its inverse counting multiplicity. Then there exist polynomials $P,Q$ and $c\in\C$, such that $f(z)=\int_{0}^{z}P(t)\exp(Q(t))dt+c$.
\end{thm}
For an entire transcendental functions $f$ with only finitely many singularities of the inverse, all of which are pre-periodic or escape exponentially, the set $P(f)$ does not accumulate in $\C$, in particular not everywhere in $\C$. Therefore if (ii) is satisfied, the function cannot be recurrent. Thus for this restricted family of functions, the question whether (i) or (ii) is true, is equivalent to the question whether $f$ is recurrent or not. As an answer to this question we get
\begin{thm}\label{class}
Let $f$ be entire and transcendental with only a finite number of singularities of its inverse, counting multiplicity, such that all these either escape exponentially or are pre-periodic, but no critical point is periodic. Then $f$ is not recurrent, if and only if all asymptotic values escape exponentially.
\end{thm}
It is remarkable that this only depends on the asymptotic values.\\
In the last chapter we discuss applications of theorem \ref{asatz} for other families, especially transcendental meromorphic functions with rational Schwarzian derivative.

\section{Basic tools}
We will use the following notation:\\
Let $f^k$ denote the $k$-th iterate, and $f^{(k)}$ the $k$-th derivative of $f$. Let ``$\meas$'' denote the Lebesgue
measure, ``$\dist$'' the Euclidean distance, and ``$\diam$'' the diameter in $\C$. Let $B(z,r)$ denote the open ball of
radius $r$ and center $z$, $B(M,\ep):=\bigcup_{z\in M}B(z,\ep)$ for $M\subset\C$ and $\D(r):=\C\Lbac B(0,r)$. For a square $S$ let $r S$ denote the square with the same center,
satisfying $\diam(rS)=r\diam(S)$.
For a conformal map $f:D\to\C$ we call $\sup_{z,w}\left|\frac{f'(z)}{f'(w)}\right|$ its \emph{distortion}.\\
We state the well known Koebe distortion theorem as it may be found in \cite{pommerenke}.
\begin{thm}[Koebe]\label{koebedist}
Suppose $f:B(0,1)\to\C$ be conformal with $f(0)=0$ and $f'(0)=1$ and $z\in B(0,1)$. Then
\begin{eqnarray}
\frac{1-|z|}{(1+|z|)^3}\le|f'(z)|  \le \frac{1+|z|}{(1-|z|)^3},
\end{eqnarray} 
\begin{eqnarray}
\frac{|z|}{(1+|z|)^2}\le|f(z)|  \le \frac{|z|}{(1-|z|)^2},
\end{eqnarray} 
\begin{eqnarray}
\frac{1-|z|}{1+|z|}\le\left|z\frac{f'(z)}{f(z)}\right|  \le \frac{1+|z|}{1-|z|}.
\end{eqnarray} 
\end{thm}
This implies in particular the following fact, which is known as Koebe's $\frac{1}{4}$-theorem.
\begin{cor}[Koebe]
Let $f$ be as before. Then
\begin{eqnarray}\label{koebe14}
B\left(0,\frac{1}{4}\right)\subset f\left(B(0,1)\right).
\end{eqnarray}\end{cor}
Much easier to show is the following property, which will be sufficient for most of our purposes. 
\begin{lem}
Let $f:B(z_0,r)\to\C$ be holomorphic. Then 
\begin{eqnarray}\label{qs3}
B\left(f(z_0),\inf_{z\in B(z_0,r)}|f'(z)|r\right)\subset f(B(z_0,r)).
\end{eqnarray} 
\end{lem}
\proof
We can assume that $f$ has no critical points.
We consider the straight path $\gamma$ from $f(z_0)$ to the closest boundary point of the image. The pre-image of $\gamma$ contains a path $\gamma'$ connecting $z_0$ with the boundary of $B(z_0,r)$, which is mapped by $f$ one to one onto $\gamma$. Since the length of $\gamma'$ is at least $r$, the length of $\gamma$ is at least $\inf_{z\in B(z_0,r)}|f'(z)|r$ .\qed
More than on disks we will be interested on the distortion on squares. From Koebe's distortion theorem one can obtain similar estimates for squares. The following lemma will be sufficient for our purpose and follows from Koebe's distortion theorem. However one could also prove this more elementary, using normal families. 
\begin{lem}  \label{Kc} For any $0<c<1$, there exists a $K_c>0$, such that for any holomorphic function, which is injective on some square $S$, the distortion of its restriction to $c S$ is bounded by $K_c$. Moreover $K_c$ tends to one, if $c$ tends to zero.
\end{lem}
The following two lemmas follow directly from the transformation formula.
\begin{lem} Suppose that the distortion of the conformal map $f$ is bounded by $K$. Let $D$ and $M$ be measurable subsets of its domain of definition, such that $\meas(D)>0$.  Then\begin{eqnarray}\frac{\meas(M\cap D) }{\meas(D)}\le\frac{ K^2 \meas(f(M)\cap f(D))}{\meas(f(D))}.\label{qs2}\end{eqnarray}
\end{lem}
The term on the left side of (\ref{qs2}) is called the \emph{density} of $M$ in $D$.
\begin{lem} Let $D$ be a $K$-quasi-square and $\ep>0$. Then \begin{eqnarray}\meas(D)\ge \frac{\diam(D)^2}{2 K^2}\;\text{and}  \;\meas(D\cap B(\partial D,\ep))\le 4 \ep K^2 \diam(D)\label{qs1} .\end{eqnarray}
\end{lem}
Finally we state a tool, which we will frequently use to obtain injectivity of a function on certain sets. It is a corollary of the so called monodromy theorem. This may be found in most function theory books as \cite{conway}.
\begin{lem}\label{injectivity}
Let $D'\subset D\subset \C$ be domains, $f:D\to\C$ holomorphic, such that all singularities of the inverse of $f$ are contained in the unbounded component of $\C\Lbac f(D')$. Then $f$ is injective on $D'$. 
\end{lem}
To avoid confusion we include a definition of a \emph{singularity} of $f^{-1}$.
\begin{defn}\label{defsing}
Let $D\subset\Ch$ be a domain, $f:D\to \Ch$ be meromomorphic and $s\in\Ch$. Then $s$ is called a \emph{singularity of $f^{-1}$} if there exist
\begin{itemize}
\item a smooth function $\gamma:[0,1]\to \Ch$ with $\gamma(1)=s$,
\item a domain $U\subset \Ch$ with $\gamma([0,1))\subset U$,
\item and a branch $\phi$ of the inverse of $f$ on $U$, i.e.: $\phi:U\to D$ meromorphic with $f(\phi(z))=z$,
\end{itemize}
such that there is no domain $V\subset \Ch$ with $\gamma([0,1])\subset V$, and no branch $\psi$ of the inverse of $f$ on $V$, that coincides with $\phi$ on that component of $U\cap V$, that contains $\gamma([0,1))$.\\
We denote the set of singularities of $f^{-1}$ with $\sing(f^{-1})$.
\end{defn}
In the literature sometimes the closure of this set is denoted by the same name. However if this set is finite, which is the case for all function which we will consider, this makes no difference. 
Studying the set $A:=\bigcap_{t\in(0,1)}\overline{\phi((t,1)) }$, one can classify these as follows. 
\begin{thm}
Let $D,f$ be as above and $s\in\sing(f^{-1}))$. Define $\gamma, U$ and $\phi$ as in definition \ref{defsing}. Then one of the following cases holds:
\begin{itemize}
\item There exists $z\in D$ with $f(z)=s$ and $\phi(\gamma(t))\to z$ as $t\to 1$. If neither $z$ nor $s$ coincides with $\infty$ it follows $f'(z)=0$. 
\item $ \dist(\phi(\gamma(t)),\partial D)\to 0$ as $t\to 1$.
\end{itemize}
In the first case $s$ is called a \emph{critical} value and in the second case an \emph{asymptotic} value.
\end{thm}
It is evident, that the pre-image of a neighborhood of an asymptotic value of an entire function must contain an unbounded component. Thus the multiplicity, as defined in the introduction, is always at least one.

\section{Non-recurrence}
\begin{figure}
\begin{center}
\scalebox{0.6}{\input{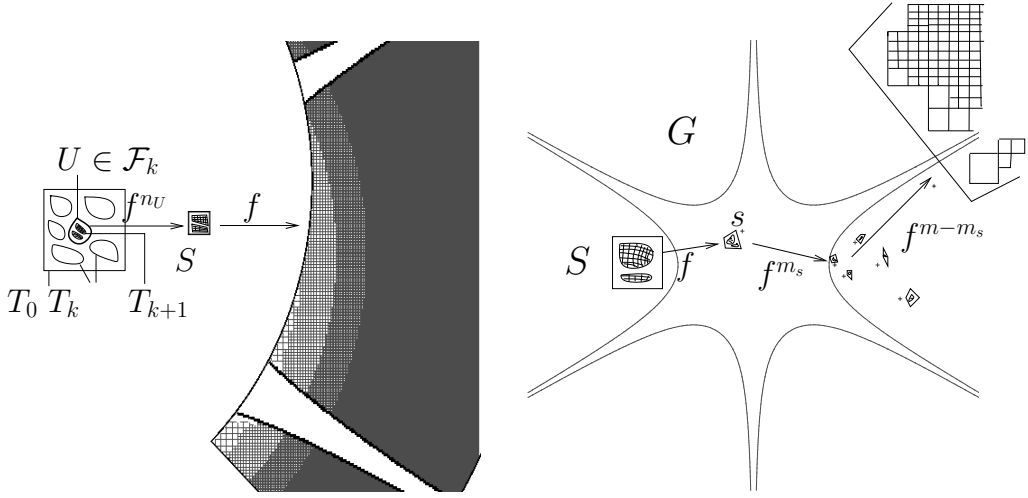}}
\end{center}
\caption{Models for the construction in both cases}\label{cases12}
\end{figure}

We follow the ideas used by M. Rees's for the exponential function. We obtain a set of points with positive measure, whose
orbits are not dense in $\C$, and therefore rule out case (i). Therefore this provides a set of sufficient
conditions for case (ii). In order to allow a wide
application, and hoping for further generalizations, we state our
theorem as general as possible. This causes a very
technical outlook.
\begin{thm}\label{asatz}
Let $f$ be meromorphic, $A\subset \C$ finite and $\G \subset \C$, such that
\begin{enumerate}
\item there exists $\ep>0$, such that the map \begin{displaymath}\overline{s}:G\to A\cup\{0\};z\mapsto\Big\{ \begin{array}{llll}s&\text{if} &  \exists s\in A:|f(z)-s|\le\exp(-|z|^{\ep})\\ 0&\text{if} &|f(z)|\ge \exp(|z|^{\ep}) \\
\end{array}
\end{displaymath}is well defined and there are $\delta_1,\delta_2\in\R$, such that for all $z\in\G$ holds \[|z|^{\delta_1}\le
\left|\frac{f'\left(z\right)}{f\left(z\right)-\overline{s}(z)}\right|\le
|z|^{\delta_2};\]
\item there exist $B>1$ and $\beta\in(-\infty,1)$, such that for every measurable $D\subset \{z: \dist\left(z,\C\Lbac\G
\right)\le 2 |z|^{-\delta_1}\}$ \[\meas(D) \le
B\diam\left(D\right)\sup_{z\in D}|z|^{\beta};\]
\item $f^m(s)\overset{m\to \infty}{\to} \infty$ and $B\left(f^m\left(s\right),2\big|f^m\left(s\right)\big|^{\tau}\right)\subset
\G $ for some $\tau> \beta$, almost all $m\in\N$ and all $s\in A$.
\end{enumerate}
Then the set $T(f):=\{z :\omega(z)\subset \overline{O^+(A)}\}$ has positive measure.\\Furthermore there exists
$M>0$, such that for any square $T_0\subset\left\{z: \dist\left(z,\C\Lbac\G \right)>|z|^{-\delta_1}\right\}$ with
$M_0:=\inf_{z\in T_0}|z|>M$ and $\diam(T_0)\ge M_0^{-\delta_2}$:
\[\frac{\meas(T(f)\cap T_0)}{\meas(T_0)}\ge 1-\exp\left(-\eta M_0^{\ep}\right),\]
where $\eta:=\frac{\tau-\beta}{\max\{1,2-2\tau\}}>0$.
\end{thm}
\begin{rmk}\label{remk} 
\rm We would like to know that $T(f)\subset J(f)$. Since the orbits of all points in $T(f)$ accumulate at infinity, the only components that could possibly intersect $T(f)$ are Baker domains and wandering domains.\\
There are various families in which these do not occur. For the family, which we consider in the next chapter, the absence of wandering has been shown by I. N. Baker in \cite{baker}. For entire functions with only finitely many singularities of the inverse this has been shown by A. Eremenko and M. Lyubich in \cite{ere} and by L. R. Goldberg and L. Keen in \cite{goldbergkeen}. For meromorphic functions with the same property this has been shown by I. N. Baker, J. Kotus and Y. L\"u in \cite{bakerkotuslu}. The absence of Baker domains has been shown for entire functions with a bounded set of singularities of the inverse by A. Eremenko and M. Lyubich in \cite{ere} and for meromorphic functions for which this set is finite by P. J. Rippon and G. M. Stallard in \cite{ripponstallard}. Moreover in \cite{bakerinf} I. N. Baker obtained an upper estimate of the growth of $|f^n(z)|$ for a point $z$ in the Baker domain of an entire function, which is not compatible with the iterated exponential escape, which we will find in the proof below for points escaping to $\infty$ in $T(f)$. Similar estimates implying the same for meromorphic functions have been obtained in \cite{bakerkotuslu3}. \\ 
It also makes sense to chose $A=\emptyset$. Then we obtain sufficient conditions for $\meas(I(f))>0$, if $I(f)$ denotes the set of escaping points (see \ref{if}).
\end{rmk}
\proofof{theorem \ref{asatz}}
 From our conditions (b) and (c) one can deduce that
$-\delta_1<\beta<\tau<1$. We note that for any $M>0$ a
sufficiently large choice of $B$ allows us to choose $G$ such that
$\G\cap B(0,M)=\emptyset$.
For all
$s\in A$ we define
\begin{eqnarray}m_s:=\max(\{m\in\N:f'(f^{m-1}(s))=0\}\cup\{0\})
\end{eqnarray}
and 
\begin{eqnarray}
k_s:=\min\{k\in\N:(f^{m_s})^{(k)}(s)\not=0\}.
\end{eqnarray}
The distortion constant $K_c$ from lemma \ref{Kc} tends to one as $c$ tends to zero. Thus for $c>0$ small enough we have $\frac{cK_c}{4}<1$.
Since $A$ is finite one can even find $c>0$, such that  $k_s\arcsin(\frac{c K}{4})<\pi$ holds for all $s\in A $ and $K:=K_c$.
Suppose that $0<\delta$ is small. In fact it turns out that $\delta <\frac{(\tau-\beta)(1-\tau)}{6-5\tau-\beta}$ is sufficient for all requirements needed. Similarly chose $M>0$ sufficiently large, satisfying many bounds appearing
throughout the proof. For now we only require the following two properties: Firstly for any $M_0>M$ the series, defined
by $M_{k+1}:=\exp\left( \min\{1,\frac{1}{2-2\tau}\} M_k^{\ep}\right)$, tends to infinity fast
enough, such that
\begin{eqnarray}\label{defm0}
\prod_{k\in\N}\left(1-\frac{1}{4}M_k^{\beta-\tau}\right)\ge
1-M_1^{\beta-\tau}.
\end{eqnarray} Secondly there are no
critical points in $A_k$ for all $k\in \N\cup\{0\}$ where
\begin{eqnarray*}
A_k&:=& \left(\D\left({\frac{1}{2}M_{k+1}^{\frac{1}{1+2\delta_2-\delta_1+3\delta}}}\right)\cap G\right)\\ &\cup& \bigcup_{0\le l\le m_s} B\left( f^l(s),M_{k+1}^{\delta-1}\right)\Lbac\{f^l(s)\}\\&\cup& \bigcup_{l>m_s}B(f^l(s), a_{k,s,l}),
\end{eqnarray*}
and $a_{k,s,l}:=\sup\{|f^j(s)|^{-\delta_2}:j\ge l, |f^{j+1}(s)|\ge M_{k+1}^{\frac{1}{1+2\delta_2-\delta_1+3\delta}}\}$. Of course, at this point we only need to study $A_0$ since the $A_k$ are descending. This does not contain any critical points for $M_0$ large enough, since those do not accumulate in $\C$ and, with
condition (a), $G$ does not contain any critical points. We note that due to (c) every $s\in A$ escapes in $\G$ exponentially fast, such that $a_{k,s,l}\le |f^l(s)|^{-\delta_2}$ for large $l$.
Now let $T_0$ and $M_0$ be as in the theorem. Let $\S$ be a family of
disjoint open squares $S\subset\{z:\dist(z,\C\Lbac G)\ge |z|^{-\delta_1}\}$ satisfying 
\begin{eqnarray}\frac{c}{8}\left(\inf_{z\in \frac{1}{c}S}|z|\right)^{-\delta_2}\le
\diam\left(S\right)\le \frac{c }{2}\left(\sup_{z\in
\frac{1}{c}S}|z|\right)^{-\delta_2},\label{diamS}\end{eqnarray} whose union cover $\{z\in
G:\dist(z,\C\Lbac G)\ge 2|z|^{-\delta_1}\}$ up to measure zero, such that $\overline{T_0}=\bigcup_{S\in X}\overline{S}$
for some finite $X\subset \S $. A typical picture could look like figure \ref{familyS}.
\begin{figure}
\begin{center}
\scalebox{0.6}{\input{familieS.pstex_t}}
\end{center}
\caption{The family $\S$ }\label{familyS}
\end{figure}
We can get this by covering the whole plane with open squares of a constant diameter,
beginning with $T_0$, cutting these into four until their diameter satisfies the upper bound, and throwing away those
intersecting $\{z:\dist(z,\C\Lbac G)\le \frac{|z|^{-\delta_1}}{2}\}$. For $M$ large enough and $G\cap
B(0,M)=\emptyset$ our squares also satisfy the lower bound. We prove the measure estimate in the theorem for all elements of $\S$ including the ones in $X$. This implies this estimate also for $T_0$. Thus we proceed with an element of $\S$, which we again call $T_0$.\\
With the estimates of condition (a) one can show that if $|z_0|$ is large enough and $B(z_0,|z_0|^{-\delta_1})\subset G$ then \begin{eqnarray}\label{inject}\text{$f$
is injective on $B(z_0,\frac{|z_0|^{-\delta_2}}{4})$.}\end{eqnarray} To see this we use lemma \ref{injectivity} and show first that
\begin{eqnarray}\label{subsub}f\left(B\left(z_0,\frac{|z_0|^{-\delta_2}}{4}\right)\right)\subset B\left(f(z_0),\frac{3|f(z_0)-\overline{s}(z_0)|}{8}\right) \subset f\left(B(z_0,|z_0|^{-\delta_1})\right).\end{eqnarray}
If the first inclusion was not true, we would find $z\in B\left(z_0,\frac{|z_0|^{-\delta_2}}{4}\right) $ with $|f(z)-f(z_0)|\ge \frac{3|f(z_0)-\overline{s}(z_0)|}{8}$. We choose $|z-z_0|$ minimal with this property, such that for $x\in(z,z_0):=\{(1-t)z+tz_0:0<t<1\}$ we have $|f(x)-\overline{s}(z_0)|\le |f(x)-f(z_0)|+|f(z_0)-\overline{s}(z_0)|\le\frac{11}{8} |f(z_0)-\overline{s}(z_0)| $. The mean value theorem provides $x\in(z,z_0)$ with\[|f'(x)|\ge\frac{|f(z)-f(z_0)|}{|z-z_0|}\ge \frac{3}{2}|f(z_0)-\overline{s}(z_0)||z_0|^{\delta_2}\ge\frac{12}{11}|f(x)-\overline{s}(z_0)||z_0|^{\delta_2}.\]
This contradicts (a), since for $z_0$ large enough $|x-z_0|\le |z_0|^{-\delta_2}$ is very small, such that $\overline{s}(x)=\overline{s}(z_0)$ and $|x|> (\frac{11}{12})^{-\delta_2}|z_0|$.\\ 
The second inclusion from (\ref{subsub}) follows from the the fact that there are no critical points in
$B(z_0,|z_0|^{-\delta_1})\subset G $. Thus we may extend the branch of $f^{-1}$, mapping $f(z_0)$ to
$z_0$, along any path in $B(f(z_0),\frac{3|f(z_0)-\overline{s}(z_0)|}{8})$ as long as the image stays in $B(z_0,|z_0|^{-\delta_1})$. As above the mean value theorem and condition (a) assures this, since for $x\in B(z_0,|z_0|^{-\delta_1})$ with $f(x)\in  B(f(z_0),\frac{3|f(z_0)-\overline{s}(z_0)|}{8})$ we know that $|f(x)-\overline{s}(z_0)|\ge \frac{5|f(z_0)-\overline{s}(z_0)|}{8} $, such that \begin{eqnarray}|f'(x)|\ge |f(x)-\overline{s}(z_0)||x|^{\delta_1}\ge\frac{1}{2}|f(z_0)-\overline{s}(z_0)||z_0|^{\delta_1}.
\end{eqnarray} This implies that the image of any path in $B(f(z_0),\frac{3|f(z_0)-\overline{s}(z_0)|}{8}))$ stays in fact inside $B(z_0,\frac{3}{4}|z_0|^{-\delta_1})$.\\ 
From this and (\ref{subsub}) follows that $f$
is injective on $B(z_0,\frac{|z_0|^{-\delta_2}}{4})$, as claimed. Together with (\ref{diamS}) this implies that the distortion of $f$ on any $S\in\S$ is
bounded by $K$.\\
Starting with $F_0:=\{T_0\}$ and $n_0(T_0):=0$, we will define for every $k\in\N$ a family
$\F_k$ of disjoint simply connected domains and functions $n_k:\F_k\to\N$, such that the sets $T_k:=\bigcup \F_k=\bigcup_{F\in\F_k}F$ form a decreasing series with the following properties for every $U\in \F_k$ and the corresponding $V\in\F_{k-1}$ with $U\subset V$:\
\begin{itemize}
\item[(i)] $\D({M_k})\supset f^{n_k\left(U\right)}\left(U\right)\in \S$ and $ \frac{1}{c}f^{n_{k}\left(U\right)}\left(U\right)\subset f^{n_{k}(U)}{(V)}$,
\item[(ii)] $f^j(V)\subset A_k$ for every $n_{k-1}(V)< j< n_{k}(U)$,
\item[(iii)] $\meas\left(V\cap \bigcup\F_{k}\right)\ge \left(1-\frac{1}{4} M_{k}^{\beta-\tau}\right)\meas(V).$
\end{itemize}
The condition (ii) implies that $\omega\left(z\right)\subset
\bigcap_{k\in\N}A_k=O^+(A)\cup\{\infty\}$ for all $z\in
T:=\bigcap_{k\in\N} T_k$. Having (iii) for each component of $T_k$, namely the elements of $\F_k$,
implies that
$\meas\left(T_{k}\right)>\left(1-\frac{1}{4}M_k^{\beta-\tau}
\right)\meas(T_{k-1})$, which, together with the exponential growth of $M_k$,
guarantees that
\[\meas\left(T\right)\ge\left( \prod_{k=1}^{\infty} \left(1-\frac{1}{4}
M_k^{\beta-\tau} \right)\right) \meas(T_0).\] Together with
(\ref{defm0}) this is the second part of our claim and thus completes the proof of theorem \ref{asatz}.\\ It remains to construct the sequences. We will do so inductively and assume the existence of appropriate $\F_k$ and $n_k$
for some $k\in\N$. We note that the starting step of the induction works the same way as any other step, such hat we do not consider it separately. Let $U\in\F_k$. Then $S:=f^{n_k(U)}(U)\in\S$. Due to condition (ii) and the
fact that there are no critical points in $A_0$, one can extend the inverse of $f^{n_k(U)}\big|U$ to $\frac{1}{c}S$ and its distortion on $S$ is bounded by
$K$. Furthermore $S\subset\D({M_k})$, such that we can consider
the following cases separately:\\
\quad\\
Case 1: $f(S)\subset\D({M_{k+1}})$. \\
We define $\F:=\{R\in \S :
\frac{1}{c}R\subset f(S)\}$, $\F_U:=\{(f^{n_k(U)+1}|U)^{-1}(R): R \in \F\}$ and $n_{k+1}(V):=n_k(U)+1 $ for all $W\in
\F_U$. See figure \ref{cases12}. Then (for $\F_U$ in place of $\F_{k+1}$) property (i) holds by definition, while property (ii) is trivial. Since
$f|S$ is injective and its distortion is bounded by $K$, $f(S)$ is a $K$-quasi-square with
\begin{eqnarray}\label{aund3}\diam(f(S))\ge \frac{1}{\sqrt{2}}\diam(S)\inf_{z\in S}|f'(z)| \overset{ }{\ge} \sup_{z\in S}|f(z)|^{1-\delta} \end{eqnarray}
for $M_0$ large enough. Here the last inequality holds since the term $|f'(z)|$ is due to (a) of magnitude $|f(z)|\ge \exp(|z|^\ep)$, which is far larger than all other factors that appear, such that those may be canceled by $|f(z)|^\delta$. Also the infimum may be substituted by the supremum, since the distortion is bounded by $K$.\\
 By definition of $\S$ and $\F$ the
set $f(S)\Lbac \bigcup\F$ is contained in the union of
$\partial\bigcup \S $, which has measure zero, and small
neighborhoods of $\partial f(S)$ and $ \C\Lbac G$. More precisely
we have
\begin{eqnarray}\label{firstterm}
\meas\left( f(S)\Lbac\bigcup\F\right)&\le& \meas\Big(\big\{z\in
f(S):\dist(z,\partial f(S) )\le |z|^{-\delta_2}\big\}\Big)\nonumber\\
&+&\meas\Big(\big\{z\in f(S):\dist(z,\C\Lbac
G)\le 2|z|^{-\delta_1})\big\}\Big).
\end{eqnarray}
With condition (b) we can control the second term on the right by
 \[\meas(\{z\in f(S):\dist(z,\C\Lbac f(S))\le 2 |z|^{-\delta_1}
\}\le B\diam(f(S))\sup_{z\in f(S)}|z|^{\beta} .\] $f(S)$ is a $K$-quasi-square. Therefore the
measure of an $r$-neighborhood of the boundary of $f(S)$ is, due to (\ref{qs1}), at most $4 rK^2\diam(f(S))$ and $\meas(f(S))\ge
\diam(f(S))^2/(2K^2)$. Since the set of the first term in (\ref{firstterm}) is contained in a  $\sup_{z\in
f(S)}|z|^{-\delta_2} $-neighborhood of $\partial f(S)$ and $-\delta_2<-\delta_1<\beta$, we obtain, using (\ref{aund3}), that \begin{eqnarray}\label{meascase1} \frac{
\meas(f(S)\Lbac \bigcup\F)}{\meas(f(S))}&\le&
\frac{2K^2(4 K^2 \sup_{z\in f(S)}|z|^{-\delta_2}+  B\sup_{z\in
f(S)}|z|^{\beta})}{\diam(f(S))}\nonumber\\
&\overset{}{\le}& 8 B K^4 \sup_{z\in f(S)}|z|^{\beta+\delta-1}\nonumber\\&\le& 8 B K^4 M_{k+1}^{\beta+\delta-1}
\end{eqnarray} for $M$ large enough. As mentioned, the distortions of $f^{n_k(U)}|U$ and $f|S $ are bounded by $K$. Therefore the distortion of $f^{n_{k+1}(V)}|U$ is bounded by $K^2$ and we get \[\frac{\meas(U\Lbac \bigcup \F_U)}{\meas(U)}\le\frac{K^4\meas(f(S)\Lbac\bigcup \F  )}{\meas(f(S))},\]
which, together with (\ref{meascase1}), implies property (iii) for $\delta$ small and $M$ large enough.\\
Case 2: $f(S)\subset B(s,M_{k+1}^{-1})$ for some $s\in A$.\\
We will study the behavior of a certain number of iterates on of $f$ on
$S$. See also figure \ref{cases12}. We begin with the first iterate. Let $w$ be the center of
$S$. For $z\in S$ (\ref{diamS}) implies $|z-w|\le \frac{c}{4}|w|^{-\delta_2}$ and (a) implies $|f'(w)|\le |f(w)-s||w|^{\delta_2}$. The mean value theorem provides $x\in
[z,w]$ with
\[|f(z)-f(w)|\le |f'(x)||z-w|\overset{}{\le} K|f'(w)||z-w|\overset{}{\le} \frac{c}{4} K |f(w)-s|.\]
Thus we know that\begin{eqnarray}\label{it11} f(S)\subset
B\left(s,(1+\frac{c K}{4})|f(w)-s|\right)\Lbac B\left(s,(1-\frac{c K}{4})|f(w)-s|
\right).
\end{eqnarray}
$S$ contains the disc $B(w,\frac{c}{16
\sqrt{2}|w|^{\delta_2}})$. Thus (\ref{qs3}) together with (a) implies that
\begin{eqnarray}\label{it12}
B\left(f(w),\frac{c|f(w)-s|}{16
\sqrt{2}K|w|^{\delta_2-\delta_1}}\right)\subset f(S).\end{eqnarray}
Next we consider those iterates, in which we cannot avoid critical
points. We do this in terms of the power series
\[f^{m_s}(z)=f^{m_s}(s)+(f^{m_s})^{(k_s)}(s)(z-s)^{k_s}+O\left((z-s)^{k_s+1}\right).\]
This provides good estimates for $f^{m_s}$ and its derivative, if $|z-s|$ is very small, which is the case for $z\in f(S)$ since $|f(w)-s|\le\exp(-|w|^{\ep})$. The only purpose of our choice of $c$ was to achieve that the diameter of $f(S)$ is small enough to ensure that $f^{m_s}(s)$ lies in the unbounded component of $\C\Lbac f^{m_s+1}(S)$. The reader may prefer to convince himself that this goal is achievable by a sufficiently small choice of $c$, instead of checking that our concrete choice above is sufficient. Thus lemma \ref{injectivity} implies that $f^{m_s}$ is injective on $f(S)$. The ratio of the outer and inner radii of the annulus in (\ref{it11}) is $C:=\left(\frac{1+\frac{c K}{4}}{1-\frac{c K}{4}}\right)$. Thus the image lies in an annulus whose ratio of those radii is very close to $C^{k_s}$ and the distortion is bounded by any constant greater than $C^{k_s-1}$, say $C^{k_s}$. Using the factor $1\pm\frac{cK}{4}$ for the error term of the power series we can deduce from (\ref{it11}) and (\ref{it12}) that 
\begin{eqnarray}\label{ms+1}
f^{m_s+1}(S)&\supset&B\big(f^{m_s+1}(w),\big|\frac{k_s(f^{m_s)^{(k_s)}}(s)
(1-\frac{c K}{4})^{k_s} (f(w)-s)^{k_s}}{16\sqrt{2}K|w|^{\delta_2-\delta_1}} \big|
\big)\end{eqnarray}
and\begin{eqnarray}\label{ms+12}
f^{m_s+1}(S)&\subset&
B\left(f^{m_s}(s),\left|(f^{m_s})^{(k_s)}(s)(1+\frac{c
K}{4})^{k_s+1}(f(w)-s)^{k_s}\right|\right).
\end{eqnarray}
Here all but the term $|f(w)-s|$ do not depend on $k$. One could get similar estimates for $1\le l\le m_s+1$, that imply $f^l(S)\subset B(f^{l-1}(s), M_{k+1}^{\delta-1})\Lbac \{f^{l-1}(s)\}$ for $M$ large enough. This implies that $f^l(S)$ is contained in $A_k$ or, more precise, in the middle term of its definition.\\ Next we consider the
maximal number of iterates, where we can assure injectivity and bounded distortion. Due to (c) the set
$B(f^m(s), 16|f^m(s)|^{-(\delta_2+\delta)})$ is contained in $ \{z: \dist(z,\C\Lbac G)\ge |z|^{-\delta_1}\}\cap
B(f^m(s), |f^m(s)|^{-\delta_2}) $ for $m$ large enough. For $m$ large $f^m(s)|^{\delta}>64$ such that, due to (\ref{inject}), $f$ restricted to this set is injective. By Koebe's $1/4$-theorem we get
\begin{eqnarray}
f\left(B\left(f^{m}(s),8|f^{m}(s)|^{-(\delta_2+\delta)}\right)\right)&\supset&
B\left(f^{m+1}(s),\frac{2|f'(f^{m})(s)|}{|f^{m}(s)|^{(\delta_2+\delta)}}\right)\nonumber\\
&\overset{}{\supset}&
B\left(f^{m+1}(s),\frac{2|f^{m+1}(s)|}{|f^{m}(s)|^{\delta_2-\delta_1+\delta}}\right).\label{mstom}
\end{eqnarray}
Here the last inclusion follows with (a) and (c).
Due to condition (c) we know that $f^m(s)$ escapes to $\infty$ in $G$. For $z\in G$ (a) implies $|f'(z)|\ge|f(z)||z|^{\delta_1}$. Thus $|f'(f^m(s))|\to\infty$ as $m\to\infty$. Thus for $m$ large enough, $m_s\le l\le m$ and $r>0$ small we know that $ f^{m-l}\text{ is expanding on  } B(f^l(s),r).$ Consequently the component of
$(f^{m-l})^{-1}\left(B\left(f^{m}(s),8|f^{m}(s)|^{-(\delta_2+\delta)}\right)\right)$, containing $f^l(s)$, is contained in
$B\left(f^{l}(s),8|f^{m}(s)|^{-(\delta_2+\delta)}\right)$, which does not contain critical points. This allows us to extend the inverse $g$ of
$f^{m-m_s+1}$, mapping $f^{m+1}(s)$ to $f^{m_s}(s)$,
to
$B\left(f^{m+1}(s),\frac{2|f^{m+1}(s)|}{|f^{m}(s)|^{\delta_2-\delta_1+\delta}}\right)$. Thus the distortion on half the ball is bounded by some constant $\tilde{K}$. One could use lemma \ref{Kc} to obtain $\tilde{K}=K_{\frac{1}{\sqrt{2}}}$ or the original distortion theorem \ref{koebedist} to obtain $K_{\frac{1}{\sqrt{2}}}=81$. In any case we get
\begin{eqnarray*}
&\dist\left(f^{m_s}(s),\partial g\left(B\left(f^{m+1}(s),\frac{|f^{m+1}(s)|}{|f^{m}(s)|^{\delta_2-\delta_1+\delta}}\right) \right)\right)\\
\ge&\displaystyle\frac{|f^{m+1}(s)|}{|f^{m}(s)|^{\delta_2-\delta_1+\delta}}\inf_{z\in B\left(f^{m+1}(s),\frac{|f^{m+1}(s)|}{|f^{m}(s)|^{\delta_2-\delta_1+\delta}}\right)}|g'(z)|\\
\overset{\text{}}{\ge}&\displaystyle \frac{|f^{m+1}(s)|}{\tilde{K}|f^{m}(s)|^{\delta_2-\delta_1+\delta}} |g'(f^{m+1}(s))|\\
\overset{\text{}}{=}&\Bigg|\displaystyle\frac{f^{m+1}(s)}{\tilde{K}f^{m}(s)^{\delta_2-\delta_1+\delta} (f^{m-m_s+1})'(f^{m_s}(s)) }\Bigg|\\
\overset{\text{}}{=}&\Bigg|\displaystyle\frac{f^{m+1}(s)}{\tilde{K}f^{m}(s)^{\delta_2-\delta_1+\delta} \prod_{i=m_s}^{m}f'(f^i(s)) }\Bigg|\\
\overset{}{\ge}&\Bigg|\displaystyle\frac{1}{\tilde{K}f^{m}(s)^{1+2\delta_2-\delta_1+\delta}f^{m_s}(s)^{\delta_2} \prod_{i=m_s+1}^{m-1}f^i(s)^{1+\delta_2} }\Bigg|\\
\overset{\text{}}{\ge}&|f^{m}(s)|^{-(1+2\delta_2-\delta_1+2\delta)}
\end{eqnarray*}
for large $m$. We define $m$ as the greatest natural number that
satisfies
\[|f^{m-1}(s)|^{-(1+2\delta_2-\delta_1+2\delta)}\ge |(f^{m_s})^{(k_s)}(s)|(1+\frac{c K}{4}
)^{k_s+1} |f(w)-s|^{k_s} .\]  We note that $m\to\infty$ as $|f(w)-s|\to 0$. Thus we can guarantee that $m$ is large by choosing $M$ large. This choice of $m$ guarantees together with (\ref{ms+1}) that $f^{m_s+1}(S)\subset g\left(B\left(f^{m}(s),\frac{|f^{m}(s)|}{|f^{m-1}(s)|^{\delta_2-\delta_1+\delta}}\right)\right) $. Thus $f^{m-m_s}$
restricted to $f^{m_s+1}(S)$ is injective, its distortion is bounded by $\tilde{K}$, and
\begin{eqnarray}\label{b12}
f^{m+1}(S)\subset
B\left(f^m(s),\frac{|f^m(s)|}{|f^{m-1}(s)|^{\delta_2-\delta_1+\delta}}\right)\subset B\left(f^m(s),\frac{|f^m(s)|}{2}\right) .\end{eqnarray}  The maximal choice of $m$
guarantees that
\begin{eqnarray}\label{vorlfm}
|f^m(s)|\ge |f(w)-s|^{\frac{-k_s}{1+2\delta_2-\delta_1+3\delta}}
\end{eqnarray} for $M_0$ large
enough. Together with the discussion below (\ref{mstom}), this implies that $f^k(S)$ is contained in $A_k$ for $m_s+1\le k\le m+1$. The exponential growth of $|f^m(s)|$ and our
estimates for $|f'|$ imply that
\begin{eqnarray}\label{expogr}|f^m(s)|^{1-\alpha}\le|(f^{m-k})'(f^k(s))|\le|f^m(s)|^{1+\alpha}\end{eqnarray}
 for any $\alpha>0$ and any natural
numbers $k<m$ with $m$ large enough. This allows us to cancel the
smaller factors, transferring (\ref{ms+1}) and (\ref{ms+12}) by bounded distortion of $f^{m-m_s}$ to
\begin{eqnarray}\label{m+1}
 B\left(f^{m+1}(w), |f^m (s)|^{1-\delta}|f(w)-s|^{k_s}\right)\subset  f^{m+1}(S)
\end{eqnarray} and
\begin{eqnarray}\label{m+12}f^{m+1}(S)\subset B\left(f^{m}(s), |f^m (s)|^{1+\delta}|f(w)-s|^{k_s}\right).
\end{eqnarray}
Again we distinguish between two cases:\\
Case 2.1 We have \begin{eqnarray}\label{case21}|f^m(s)|\ge |f(w)-s|^{\frac{- k_s}{1-\tau+\delta}}.\end{eqnarray}
This is stronger than (\ref{vorlfm}) and with (a) we get
\begin{eqnarray}\label{endfm} |f^m(s)|\ge
\exp\left(\frac{k_s}{1-\tau+\delta}M_k^{\ep}\right),\end{eqnarray} which together
with (\ref{b12}) implies that $f^{m+1}(S)\subset \D({M_{k+1}})$ for
$\delta$ sufficiently small. We define $\F:=
\{R\in\S:\frac{1}{c}R\subset f^{m+1}(S)\}$, \[\F_U:= \{(f^{n_k(U)+m+1}|U)^{-1}(R): R \in \F\}\;\text{and}\;n_{k+1}(V):=n_k(U)+m+1 \;\text{ for all }\;V\in \F_U.\]
Again properties (i) and (ii) follow by definition. The distortion of $f|S $ is bounded by $K$, while the distortion of $f^{m_s}|f(S)$
is bounded by $C^{k_s} $, and the
distortion of $f^{m-m_s}|f^{m_s+1}(S)$ is bounded by $\tilde{K}$.
Thus the distortion of $f^{m+1}|S$ is bounded by
$\tilde{K}KC^{k_s} $. Therefore
$f^{m+1}(S)$ is a $\tilde{K}KC^{k_s} $-quasi-square. For $M$ large enough (\ref{case21}) together
with (\ref{m+1}) imply that
\begin{eqnarray}\label{diam21}\diam(f^{m+1}(S))\ge |f^m(s)|^{\tau-2\delta} \ge \sup_{z\in
f^{m+1}(S)}|z|^{\tau-3\delta}.\end{eqnarray} As in case 1 we find that
$\meas(f^{m+1}(S)\Lbac \bigcup \F)$ is bounded above by the measure of the set
$\{z:\dist(z,\C\Lbac G)\le 2|z|^{-\delta_1}\}$, which is, due to condition (b), at most $B\diam(f^{m+1}(S))\sup_{z\in
f^{m+1}(S)}|z|^{\beta}$ plus the measure of the set $\{z\in f^{m+1}(S):\dist(z,\partial f^{m+1}(S))\le |z|^{-\delta_2} \}$, which is, with (\ref{qs1}), bounded above by $4 \tilde{K}^2K^2C^{2k_s} \diam(f^{m+1}(S))\sup_{f^{m+1}(S)}|z|^{-\delta_2}  $. Again with (\ref{qs1}) we know that $\meas(f^{m+1}(S))\ge
\diam(f^{m+1}(S))^2/(2\tilde{K}^2K^2C^{2k_s} )$.  Using $-\delta_2< \beta$ and (\ref{diam21}) we can deduce from the above that
\begin{eqnarray}\label{measc21}
\frac{\meas\left(f^{m+1}(S)\Lbac
\bigcup\F\right)}{\meas\left(f^{m+1}(S)\right)}\le 5 B
\tilde{K}^2K^2C^{2k_s} \!\!\!\!\!\!\!\!\!\!\sup_{z\in
f^{m+1}(S)}|z|^{\beta+3\delta-\tau}\le\frac{
M_{k+1}^{\beta-\tau}}{4 \tilde{K}^2K^4C^{2k_s}},
\end{eqnarray} where the last inequality holds due to (\ref{b12}) and (\ref{endfm}) for $M$ large and $\delta$ small enough. The distortion $f^{n_k(U)}|U$ is bounded by $K$. Thus (\ref{measc21}) together with (\ref{qs2}) implies         (iii).  \\
Case 2.2 $|f^m(s)|<  |f(w)-s|^{\frac{-k_s}{1-\tau+\delta}}$\\
With (\ref{m+12}) we get
\begin{eqnarray}\label{20a}f^{m+1}(S)\subset B(f^m(s),|f^m(s)|^{1+\delta} |f(w)-s|^{k_s})\subset B(f^m(s),|f^{m}(s)|^{\tau}),\end{eqnarray} which, because of condition (c), is contained in $\{z:\dist(z,\C\Lbac G)\ge|z|^{-\delta_1}\}$. We distinguish between two more cases:\\
Case 2.2.1 $\diam(f^{m+1}(S))< \frac{c}{4}|f^{m}(s)|^{-\delta_2} $\\
Then due to (\ref{inject}) $f|f^{m+1}(S)$ is injective and its distortion is bounded by
$K$. We define
\[\F_U:= \{(f^{n_k(U)+m+2}|U )^{-1}(T):T\in\S;\frac{1}{c}T\subset
f^{m+2}(S)\}\]
and $n_{k+1}(V):=n_k(U)+m+2$.Then the
property (i) and (ii) are again satisfied by definition. The bounds of the distortion of
$f|S$,$f^{m_s}|f(S)$ and $f^{m -m_s}|f^{m_s+1}(S)$ are as above.
Then $f^{m+2}(S)$ is a $\tilde{K}K^2C^{k_s} $-quasi-square and, with (\ref{m+1}), it follows that
\begin{eqnarray}\diam(f^{m+2}(S))\overset{}{\ge} |f^m(s)|^{1-\delta}|f(w)-s|^{k_s}\!\!\!\!\!\!\!\!\!\inf_{z\in f^{m+1}(S)}\!\!\!\!\!\!\!\!\!|f'(z)|\ge\!\!\!\!\!\!\!\! \sup_{z\in f^{m+2}(S)}\!\!\!\!\!\!|z|^{1-\delta}.
\end{eqnarray} Here the last inequality follows for $M$ large enough, since for $z\in f^{m+1}(S)$ the magnitude of  $|f(z)|$, and, with condition (a), also that of $|f'(z)|$, is $\exp(|f^m(s)|^{\ep})$, which,  with (\ref{vorlfm}), is far larger than the other factors. With condition (b) and (\ref{qs1}) we get as before
\begin{eqnarray}\label{meas221}
&&\frac{\meas(f^{m+2}(S)\Lbac \bigcup_{R\in \S,\frac{1}{c}R\subset
f^{m+2}(S)}R)}{\meas(f^{m+2}(S))}\nonumber\\
&\le&5B
\tilde{K}^2K^4C^{2k_s}\!\!\!\!\!\!\!\!\sup_{z\in f^{m+2}(S)}\!\!\!\!\!|z| ^{\beta+\delta-1}\nonumber\\
&\overset{}{\le}&\exp\left(\left(\frac{1}{2}\exp\left(\frac{
M_{k}^{\ep}}{1+2\delta_2-\delta_1+3\delta}
\right)\right)^\ep\right)^{\beta+2\delta-1}.
\end{eqnarray}
Here the last inequality follows with (\ref{b12}) and (\ref{vorlfm}).
Again with the distortion estimates from above and (\ref{qs2}) this implies that
\begin{eqnarray*}\frac{\meas(U\Lbac \F_U)}{\meas(U)}\le\frac{\tilde{K}^2K^{6}C^{2k_s}\meas(f^{m+2}(S)\Lbac \bigcup_{T\in \S,\frac{1}{c}T\subset f^{m+2}(S)}T)}{\meas(f^{m+2}(S))}. \end{eqnarray*}
Together with estimate (\ref{meas221}) this is far stronger than condition (iii).\\
Case 2.2.2 $\diam(f^{m+1}(S)) \ge \frac{c}{4}|f^{m}(s)|^{-\delta_2}$\\
We consider a family $\F$ of disjoint open squares $R\subset f^{m+1}(S)$ with  diameter $\exp(-|f^{m}(s)|^{\frac{\ep}{2}}) $, such that we cover all of $f^{m+1}(S) $ except a set of measure zero and a $\exp(-|f^{m}(s)|^{\frac{\ep}{2}}) $-neighborhood of the
boundary. Since $f^{m+1}(S)$ is a $\tilde{K}KC^{k_s}$-quasi-square, (\ref{qs1}) implies that
\begin{eqnarray}\label{denofF}\meas(f^{m+1}(S)\Lbac \bigcup\F)&\le& 4\tilde{K}^2K^2C^{2k_s}
\exp(-|f^m(s)|^{\frac{\ep}{2}})\diam(f^{m+1}(S)) \nonumber\\
&\overset{}{\le}& \exp\left(-\left|\frac{f^{m}(s)}{2}\right|^{\frac{\ep}{2}}\right)\end{eqnarray} 
for $M$ large enough, since, due to (\ref{20a}), again one factor, namely $\exp(|f^m(s)|^{\frac{\ep}{2}})$, is dominating all others. We define
\[\F_U:= \left\{\left((f^{n_k(U)+m+1}|U )^{-1}\circ (f|R)^{-1
}\right)(Q):R\in\F\;\text{and}\;  Q \in \S\;\text{with}\frac{1}{c}Q\subset f(R)\right\}\]
and $n_{k+1}(V):=n_k(U)+m+2$ for all $V\in\F_U$. Again properties (i) and (ii) follow
directly. The diameter of all $R\in\F$ is very small such that, due to (\ref{inject}), $f|R$ is injective and its distortion is close to one, say
bounded by $K$. With the mean value theorem we can deduce that \[\diam(f(R))\ge \frac{1}{\sqrt{2}}\inf_{z\in
R}|f'(z)| \exp(-|f^m(s)|^{\frac{\ep}{2}})\overset{}{\ge}\sup_{z\in
f(R)}|z|^{1-\delta}\] for $M_0$ large enough. Note that with (\ref{b12}) we have $|f(z)|\ge \exp(|\frac{f^m(s)}{2}|^\ep)$ for $z\in R$. With the same
arguments as above (b) and (\ref{qs1}) imply that
\begin{eqnarray*}
\frac{\meas(f(R) \Lbac \bigcup_{Q\in \S, \frac{1}{c}Q\subset
f(R)}Q)}{f(R)}&\le&5BK^2\sup_{z\in
f(R)}|z|^{\beta+2\delta-1}\\&\overset{}{\le}&5BK^2\exp\left(\beta+\delta-1\left|\frac{f^m(s)}{2}\right|^{\ep}\right),
\end{eqnarray*}
where last inequality may be deduced from \text{(a)} and (\ref{b12}).
Since the distortion of $f|R$ is bounded by $K$, we can transfer this with help of (\ref{qs2}) to $R$ loosing only a factor $K^2$. This estimate for the density  in every $R\in \F$ implies the same for their union $\bigcup\F$, which is contained in $f^{m+1}(S)$. More precisely we know that
\begin{eqnarray}\label{densinR}
\!\!\!\frac{\meas(\bigcup\F\Lbac\bigcup_{R\in \F,Q\in S, \frac{1}{c}Q\subset f(R)}(f|R)^{-1}(Q) )
}{\meas(f^{m+1}(S))}&\!\!\!\le&\!\!\!5BK^4\exp\left(\beta+\delta-1\left|\frac{f^{m}(s)}{2}\right| ^{\ep}\right).
\end{eqnarray} With the distortion estimates above with (\ref{qs2}) it follows from (\ref{denofF}) and (\ref{densinR}) that
\begin{eqnarray*}
&&\frac{\meas(U\Lbac\bigcup \F_U)}{\meas(U)}\\
&\le&\frac{\tilde{K}^2K^4C^{2k_s}\meas(f^{m+1}(S) \Lbac \bigcup_{R\in \F,Q\in S, \frac{1}{c}Q\subset f(R)}(f|R)^{-1}(Q) )}{\meas(f^{m+1}(S))}\\
&\le&\frac{\tilde{K}^2K^{4}C^{2k_s} (\meas(f^{m+1}(S)\Lbac\! \bigcup\! \F)\!+\!   \meas(\bigcup\F\Lbac\bigcup_{R\in \F,Q\in S, \frac{1}{c}Q\subset f(R)}(f|R)^{-1}(Q) )  
)}{\meas(f^{m+1}(S))}\\
&\overset{}{\le}& \tilde{K}^2K^4C^{2k_s}   \exp\left(-\left|\frac{f^m(s)}{2}\right|^{\frac{\ep}{2}}\right) \\
&+ & 5 K^8 \tilde{K}^2C^{2k_s}
\exp\left(\beta+\delta-1\left|\frac{f^{m}(s)}{2}\right| ^{\ep}\right).
\end{eqnarray*}
Together with (\ref{vorlfm}) this is again far stronger than condition (iii) for $\delta$ small and $M_0$ large enough.\\
\quad\\
The definition $\F_{n+1}:=\bigcup_{U\in\F_k}\F_U$ completes the recursive definition, such that all required properties are satisfied.\\
This completes the
proof of the theorem. \qed

\section{Entire functions}
\begin{figure}
\begin{center}
\scalebox{0.6}{\input{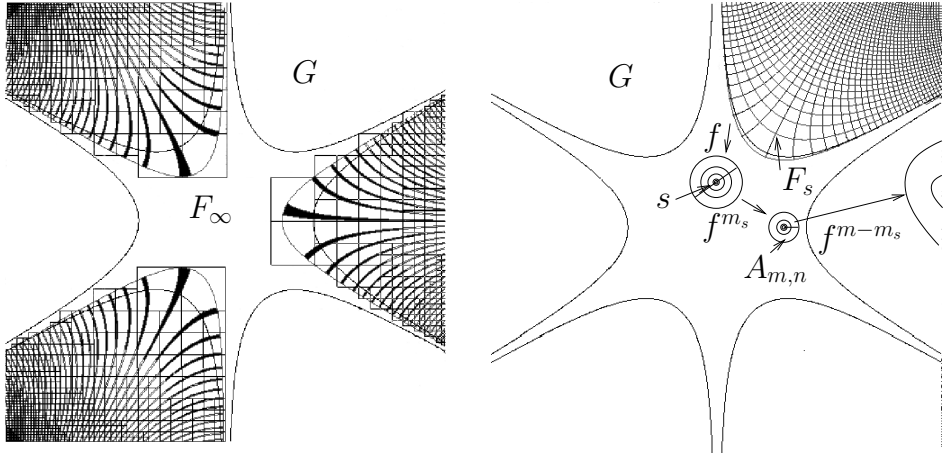}}
\end{center}
\caption{Construction of the family $\F$ }\label{familyF}
\end{figure}

In this section we will only work with functions of the same type as in theorem \ref{gsatz}. First of all we will prove some general properties and introduce some notations which will frequently occur.
In the entire chapter let $P$ and Q be polynomials with $P$ not zero and $Q$ not constant, $c\in \C$ and \begin{eqnarray}\label{gaussfunctions}f(z):=\int_{0}^z
P(t)\exp(Q(t))dt+c.\end{eqnarray}
For $k\in\{1,..,\deg(Q)\}$ define $\phi_k:=\frac{(2k+1)\pi-\arg(q)}{\deg(Q)}$, where $q$ denotes the leading coefficient of $Q(z)=q z^{\deg(Q)}+...$. For $R\to \infty $ the modulus of $\exp(Q(R\exp(\phi i))$  decreases very fast, such that $f(R\exp(\phi i))$
converges  to the point \begin{eqnarray}\label{sdach} s_k:=\lim_{R\to\infty}\int_0^{R\exp(i\phi)}
P(t)\exp(Q(t))dt +c,\end{eqnarray} which therefore is an asymptotic value of $f$. Let $A$ denote the set of asymptotic values. For $z\in\C$  choose $k$ such that $\phi_k$ is closest possible to $\arg(z)$ and define $\overline{s}(z)=s_k$. 
\begin{lem}\begin{eqnarray}
f(z)&=&\overline{s}(z)+\frac{P(z)\exp(Q(z))}{Q'(z)}+\OO(|z|^{\deg(P)-\deg(Q)})\exp(Q(z))\;\mbox{for $z\to\infty$}\label{estf}.
\end{eqnarray}\label{properties}
\end{lem}
\proof
Let $z\in \C$. We define $w:=2|z|\exp(\phi_k i)$ with the $k$ from above. Instead of integrating from $0$ to
$z$ on a straight path, one might as well go from $0$ to infinity in the direction $\phi_k$, come back the same way up to
$w$, and finally move forward to $z$. If $z,w$ are no zeros of $Q'$, one can find a path from $w$ to $z$ avoiding these
zeros, such that with integration by parts it follows that
\begin{eqnarray}
f(z)&=&s_k+ \frac{P(z)\exp(Q(z))}{Q'(z)}-\frac{P(w)\exp(Q(w))}{Q'(w)}\nonumber\\
&+&\int_{w}^{z}\frac{P'Q'-P Q''}{(Q')^2}(t)\exp(Q(t)) dt -\int_{2|z|}^{\infty}P(t\exp(\phi_k i))\exp(Q(\exp(\phi_k i))) dt.
\nonumber
\end{eqnarray}
It is easy to obtain estimates of the last three terms that imply (\ref{estf}).\qed
\begin{lem}\label{finj} Let $\delta,\delta'>0$ and $M$ large enough. Then for $z\in G:=\{z:|\Re(Q(z))|\ge |z|^\delta\}\cap D(M)$ the restriction of $f$ to $B\left(z,\frac{(1- \delta')\pi}{|Q'(z)|}\right)$ is injective. 
\end{lem} 
\proof
We use lemma \ref{injectivity}. Assume the existence of $z,w\in\G$ with $f(z)=f(w)$ and
$|z-w|<\frac{(2-2\delta')\pi}{|Q'(z)|}$. Then $f([z,w])$ is a closed curve with a singularity of $f^{-1}$ in a bounded
component of its complement. The condition $|f(z)-s|=|f(w)-s|$ together with (\ref{estf}) implies $1-\frac{\delta'}{8}<\frac{|\mbox{Re}(Q(w))|}{|\mbox{Re}(Q(z))|}<1+\frac{\delta'}{8}$ for $|w|$ large enough. Then
$1-\frac{\delta'}{4}<\frac{|\mbox{Re}(Q(x))|}{|\mbox{Re}(Q(z))|}<1+\frac{\delta'}{4}$ for every $x\in[z,w]$. Again with (\ref{estf})
we have $1-\frac{\delta'}{2}<\frac{|f(x)|}{|f(z)|}<1+\frac{\delta'}{2}$ for $|z|$ large enough. Therefore the length of the
curve $f([z,w])$ must be at least $\pi(2-\delta')|f(z)|$. On the other hand the mean value theorem does not allow this
length to exceed $|z-w|\max_{x\in[z,w]}|f'(x)| $, which is smaller than
$\frac{\pi(2-\delta'-\delta'^2)\max_{x\in[z,w]}|Q'(x)|}{|Q'(z)|}|f(z)|$ contradicting the estimate from above for
$|z|$ large enough. Thus the claim follows with lemma \ref{injectivity}.
\qed
Now we prove the first theorem of the introduction.
\proofof{theorem \ref{gsatz}} We verify the properties of theorem \ref{asatz}. Since they escape exponentially for every $s\in A$, there exists $\delta_s>0$ , such that $|f^{n+1}(s)|\ge
\exp(|f^n(s)|^{\delta_s} )$ for almost every $n\in \N$. Suppose $0<\delta<\min_{s\in A}\delta_s$. With (\ref{estf}) we have an estimate for $|f|$. 
For $\ep < \delta $ and $\delta_1<\deg(Q)-1<\delta_2$ property (a) follows, if we redefine $\overline{s}$ as zero on the part of $G$ where $\mbox{Re}(Q(z))>0$.\\
Far away from the origin $\C\Lbac G$ consists of neighborhoods around the pre-images under $Q$ of the imaginary axis, whose widths at a distance $R$ from the origin are of magnitude $R^{-\deg(Q)+1+\delta}$. With the width at a distance $R$ from the origin we mean the diameter of the largest disc that is contained in the set and whose center has modulus $R$. For $-\deg(Q)+1+\delta <\beta<1$ and $B$ sufficiently large (b) follows.\\
As mentioned above we have $|f^{n+1}(s)|\ge\exp(|f^n(s)|^{\delta_s})$ except for a finite
number of $n\in\N$. Thus the real part of $Q(f^n(s))$ is at least of magnitude
$|Q(f^n(s))|^{\frac{\delta_s}{\deg(Q)}}$.  The magnitude of the distance of $f^n(s)$ to $\C\Lbac G$ is therefore no less than 
$|f^n(s)|^{-\deg(Q)+1+\delta_s}$. If we choose $\tau$, such that $\beta<\tau < -\deg(Q)+1+ \min_{s\in A}\delta_s$, then (c) is satisfied. Now we
can apply theorem \ref{asatz} and get $\meas(T(f))>0$. Case (ii) of theorem \ref{bockclass} follows. As explained in remark \ref{remk} we know that $T(f)\subset J(f)$ due to a result of I. N. Baker \cite{baker}.\\ Assume now $\deg(Q)\ge 3$.
Then for $\delta$ small enough we have $\deg(Q)-1-\delta >1$. This implies that $\meas(\C\Lbac G)<\infty$ and if $\delta_1>1$ even that $\meas(\{z:\dist(z,\C\Lbac G)\le |z|^{-\delta_1}\})<\infty$. This follows since this set is contained in the set
\begin{eqnarray}\nonumber
B(0,M)\cup\bigcup_{k=1}^{\deg(Q)}\left\{z:\left|\arg(z)-\frac{(4k+1)\pi-2\arg(q)}{2\deg(Q)}\right|\le 2 q |z|^{\delta-\deg(Q)}\right\}.
\end{eqnarray} 
Using the transformation formula for polar coordinates, the measure of the last set is bounded by
\begin{eqnarray}\nonumber
\deg(Q)\int_M^\infty R^{1-\deg(Q)+\delta}d 2 q R= \frac{2q\deg(Q)M^{2-\deg(Q)+\delta}}{2-\deg(Q)+\delta}
\end{eqnarray}
We cover the set $\{z:\dist(z,\C \Lbac G)\ge
2|z|^{-\delta_1}\}$ with a family $\S$ of squares $S\subset \{z:\dist(z,\C \Lbac G)\ge |z|^{-\delta_1}\}$ with
$\sup_{z\in S}|z|^{-\delta_2}\le\diam(S) \le4\sup_{z\in S}|z|^{-\delta_2}$. The density of $\F(f)$ in any $S\in\S$ is,
due to theorem \ref{asatz}, at most $\exp\left(-\eta\inf_{z\in S}|z|^{\ep}\right)$.
Let $R_k :=M+k$ for all $k\in\N\cup\{0\}.$ Then it follows that
\begin{eqnarray*} \meas(F(f))&\le&
\meas\left(\{z:\dist(z,\C\Lbac G) \le
2|z|^{-\delta_1}\}\right)+\sum_{S\in\S}\meas(\F(f)\cap S)\\
&\le& \pi M^2+\frac{2q\deg(Q)M^{2-\deg(Q)+\delta}}{2-\deg(Q)+\delta}+\sum_{k\in\N}\!\!\!\!\!\!\!\!\!\!\!\!\!\sum_{\begin{array}{c}S\in\S\\R_k< \inf_{z\in S}|z|\le R_{k+1}\end{array}}\!\!\!\!\!\!\!\!\!\!\!\!\!\!\!\meas(\F(f)\cap S)\\
&\le& \pi M^2+\frac{2q\deg(Q)M^{2-\deg(Q)+\delta}}{2-\deg(Q)+\delta}+\sum_{k=1}^\infty 2\pi
(R_{k}+1)\exp\left(-\eta R_{k}^{\ep} \right)
,\end{eqnarray*}
which is finite. This gives the second part of theorem \ref{gsatz}.\qed 
This allows us to construct examples of functions $f$ with $
0<\meas(\F(f))<\infty$. For example we can arrange the parameters $a,b$ ( e.g. $a=(\frac{27\pi^2}{16})^{1/3}$, $b=\log(\sqrt{a/3}) $), such
that both critical points of the function
\[f(z)=\exp(z^3+az+b)\]
are fixed, and the only asymptotic value $0$ escapes to infinity on the real axis. Thus the Fatou set of $f$, which
consists of those two super attractive basins, has finite measure. In the figure below, where the part $\{z: |\mbox{Re}(z)|\le 2,|\mbox{Im}(z)|\le 2\}$ of the plane is displayed, these two super-attractive basins are painted dark. 
\begin{figure}[!ht!]
\begin{center}
\includegraphics[width=50mm]{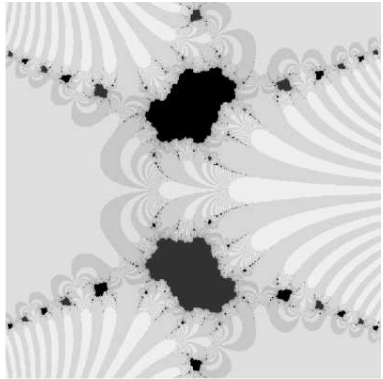}
\caption{The Fatou set of $f(z)=\exp(z^3+az+b)$ with $a=(\frac{27\pi^2}{16})^{1/3}$ and $b=\log(\sqrt{a/3})$}
\end{center}
\end{figure}
We should note that the existence of such examples is not very surprising after the construction of examples with a positive measure Julia set by C. McMullen in \cite{mcmullen}.
Also the idea of using concrete measure estimates like the one in \ref{asatz} in order to show finiteness of subsets of the Fatou set has been used before by H. Schubert, who proved in \cite{schubert}, that the measure of the Fatou set of the sine function in the strip $\{z: 0\le\mbox{Re}(z)\le 2\pi\}$ is finite, as conjectured by J. Milnor in \cite{milnor}.\quad\\
In order to prove theorem \ref{andereinclusion} we need the following lemma.
\begin{lem}\label{dengamma}
For $a\in A\cup\{\infty\}$ let $G_a$ be that part of $G$ being mapped close to $a$. Let $\Gamma\subset \C:=\bigcup_{n\in\N}(f^n)^{-1}(B(a,\ep))$ for some $a\in A$ and $\ep>0$. Then there exist positive constants $c,C$, and a family $\F$ of disjoint domains $D$, such that
\begin{eqnarray}\label{claim41}\diam(D)\le\frac{C}{\sup_{z\in D}|Q'(z)|}\;,\;\frac{\meas(\Gamma\cap D)}{\meas(D)}\ge c\; \text{and} \;\meas\left(G_s\Lbac\bigcup\F\right)=0, \end{eqnarray}
if $s=\infty$ or if $s$ is an asymptotic value that escapes exponentially. 
\end{lem}
\begin{rmk}
\rm For $M$ large and any $z\in G$ there is exactly one $a\in A$ for which $|f(z)-a|<1$ or $|f(z)|\ge \exp(\frac{M}{2})$ in which case we regard $f(z)$ as close to $\infty$. Thus the $G_a$ are well defined.   
\end{rmk}
\proof Since $a$ is an asymptotic value, we have $\lim_{R\to\infty}f(R\exp(\phi_a i))=a$ for some some
$k\in\{0,1,..\deg Q\}$ and $\phi_a=\frac{(2k+1)\pi -q}{\deg Q}$ . From (\ref{estf}) it follows that for any
$0<\delta' $ there exists $M>0$, such that \begin{eqnarray*}\left\{z:
\phi_a-\frac{(1-\delta')\pi}{2\deg(Q)}<\arg(z)<\phi_a+\frac{(1-\delta')\pi}{2\deg(Q)} )\right\}\cap D(M) \subset
f^{-1}\left(B(s,\ep)\right).\end{eqnarray*} Since (\ref{estf}) also gives good estimates for the argument of $f$ in $G_\infty$ it follows that the set
\begin{eqnarray}\label{arg}\left\{z:\phi_a-\frac{(1-2\delta')\pi}{2\deg(Q)} <\arg\left(\frac{P(z)}{Q'(z)}\right)+\mbox{Im}(Q(z))<\phi_a+\frac{(1-2\delta')\pi}{2\deg(Q)} \right\} \cap G_\infty \end{eqnarray}
is contained in the $f^{-2}\left(B(s,\ep)\right)$.
 Every component of this set is an unbounded region, whose width at the distance $R$ from the origin is at least $\frac{1-3\delta'}{\deg(Q)} 2\pi |q|^{-1}R^{1-\deg (Q)} $ for sufficiently large $R$. We refer to this regions as the ``channels'' and include the left side of figure \ref{familyF}, in which these channels are colored black,  to give an idea of their outlook. In order to be able to display the structure we had to magnify their diameter relatively to $M$. The ``gaps'' in between these channels have a width of at most $\left(2-\frac{(1-3\delta')}{\deg(Q)} \right)\pi |q|^{-1}R^{1-\deg(Q)} $, still assuming that $M$ is large. The complement of $\Gamma$ in $G_\infty$ must lie in the gaps between these channels.
 For a sufficiently small choice of $\delta'$ simple geometric arguments give for any constant $C'>\sqrt{2}\pi\left(2-\frac{1}{\deg(Q)} \right) $ a constant $c'>0$, such that for any square $S$ intersecting $G_\infty$ with $\diam(S)\ge C'(\inf_{z\in S}|Q(z)|)^{-1} $, the density of $\Gamma$ in $S$ is bounded below by $c'$. (For a sufficiently large choice of $C'$ one can choose $c'$ arbitrarily close to $\frac{1}{2\deg(Q)}$). We cover $G_\infty $ up to measure zero by a family $\F_\infty$ of squares $S$ with $C'(\inf_{z\in S}|Q'(z)|)^{-1}< \diam(S)<4C'(\sup_{z\in S}|Q'(z)|)^{-1} $. We obtain this family in a similar way as the family $\S$ in the proof of theorem \ref{asatz}: We begin with a grid of open squares covering the whole plane with a constant diameter, subdivide these into four parts until they satisfy the upper bound, and finally throw away those not intersecting $G_\infty$. Then our conditions are satisfied with $c=c'$ and $C=4C'$. The family $\F_{\infty}$ could look similar as displayed on the left side of figure \ref{familyF}. \\
Now we need to find a covering $\F_s$ of $G_s$ for every asymptotic value $s$ that escapes exponentially. We define 
$m_s\in\N\cup\{0\}$ as before minimal such that for $m\ge m_s$ the point $f^m(s)$ is no
critical point of $f$ and choose $m_0\ge m_s$ such that for $m\ge m_0$ it is none of $Q$. Then we choose $n_m$ to be the smallest natural number, for which there exists $l_m\le 4$ with
\[\left|\frac{(f^{m-m_s})'(f^{m_s}(s)) Q'(f^{m}(s))}{(f^{m-m_s-1})'(f^{m_s}(s))Q'(f^{m-1}(s))}\right|= l_m^{n_m} .\] The exponential escape
of $s$ implies that $l_m^{n_m}$ is of magnitude $|f^m(s)|$, such that $l_m\overset{m\to\infty}\longrightarrow 4$ and $n_m\overset{m\to\infty}{\longrightarrow}\infty$. For
$0\le n\le n_m$ we define \begin{eqnarray}\label{Rmn}R_{m,n}&:=& |\frac{(1-\delta')l_m^n\pi}{2(f^{m-m_s})'(f^{m_s}(s))Q'(f^m(s))}|\end{eqnarray} and for $n\not=0$ we
consider the slit annulus
\[A_{m,n}:=\{z:R_{m,n-1} <|z-f^{m_s}(s)|< R_{m,n},\;z-f^{m_s}(s)\not\in\R_{>0}\}.\] Let $\F_s$ be the family of all connected components of $(f^{m_s+1})^{-1}(A_{m,n})$, intersecting $\G_s$, for all $m\ge m_0$ and $n\in\{1,...,n_m\}$. We tried to give an idea of $\F_s$ in figure \ref{familyF}. For $M$ large enough
these components cover $G_s$ up to measure zero. Next we will verify the diameter condition.
For large $m$ the annulus $A_{m,k}$ is very close to $f^{m_s}(s)$, such that the
power-series of $f^{m_s}$ gives good estimates. If $k_s$ is the multiplicity of $f^{m_s}$ in $s$, there are $k_s$ pre-images $A'$ of $ A_{m,k}$ under $f^{m_s}$ , which
are contained in $B\left(s,r\right)\Lbac
B\left(s,r¯\right)$ with $r^{(-)}:=\frac{1+(-)\delta'}{|(f^{m_s})^{(k_s)}(s)|}\left|R_{m,n(-1)}\right|^{\frac{1}{k_s}} $. Since the ratio of the outer and inner radii of this annulus is $\frac{1+\delta'}{1-\delta'}l_m^{\frac{1}{k_s}}$ the distortion of $f^{m_s}$ on these $A'$ is bounded by any constant $C_2$, which is larger that this to the power of $k_s-1$ for $m_0$ large enough.\\
Again with estimate (\ref{estf}) one can show that
any connected component $D$ of $f^{-1}(A')$ intersecting $G_s$, which is therefore an element of $\F_s$, is a simply connected
domain with a diameter of at most $C \inf_{z\in D}\frac{1}{|Q'(z)|}$ for any $C>4^{2/k_s}+\frac{8\pi}{ k_s}$
and $M$ large enough. This follows since we can connect any two points in $A'$ with a path in $A'$ whose length is at most $(\frac{2\pi}{k_s}+l_m^{\frac{1}{k_s}})r$ and for $z\in D$ we know due to (\ref{estf}) that $|f'(z)|$ is bounded by $|Q'(z)|r^-$. Thus the diameter condition from (\ref{claim41}) is satisfied. An analogous upper estimate for $|f'|$ on $D$ implies that distortion of $f$ on such a domain $D$
is bounded by any constant $
C_3>4^{\frac{1}{k_s}}$, if $M$ is large enough. \\
It remains to show that the density of $\Gamma$ in these pre-images is again bounded away from zero by some $c>0$.\\
The diameter of $A_{m,n_m}$ is chosen in such a way that for $m$ large enough $f^{m-m_s}$ is injective on this set and its distortion is bounded by $K_{\tilde{c}}$ from lemma \ref{Kc} with $\tilde{c}<\frac{1}{\sqrt{2}l_m^{n_m-n}}$.
This follows since for $m$ large enough $f^{m-m_s-1}$ is injective on $B(f^{m_s}, R_{m-2,0})$, such that its distortion on $A_{m,n}$ is bounded by any constant larger than one, say $\frac{1-\frac{\delta'}{2}}{1-\delta'}$. Thus we have $f^{m-m_s-1}(A_{m,n})\subset B(f^{m-1}(s), \frac{(1-\frac{\delta'}{2})\pi}{2|Q'(f^{m-1}(s))|} )$. With lemma \ref{finj} we know that $f$ is injective on the twice this ball, such that the distortion estimate for $f^{m-m_s}|A_{m,n}$ follows with lemma \ref{Kc}. \\
Since $s$ escapes exponentially we may assume $\delta$ to be small enough to ensure $|f^{n+1}(s)|\ge \exp(|f^n(s)|^{2\delta})$ for large $n$, such that in particular $f^n(s)\in \G_\infty$.  \\
To show that the density of $\Gamma$
in the set $f^{m-m_s}(A_{m,n})$ is bounded below by some constant
$c>0$, we distinguish between large and small $n$.\\
If $n$ is small enough, such that $l_m^n\le |f^m(s) |^{\delta} $ is satisfied, it follows with the definition of $n_m$ that $n_m-n$ is large. Assuming $\delta<1/2$ we get $n_m-n\ge n_m/2$ for $m$ large enough. Thus the distortion of $f^{m-m_s}$ on $A_{m,n}$ is bounded by some $K_{m} $, which tends to $1$ as $m\to\infty$. Thus for $m$ large the set $f^{m-m_s}(A_{m,n})$ is close to the annulus $B\left(f^{m}(s),
\left|\frac{(1-\delta')l_m^n\pi}{2Q'(f^m(s))}\right|\right)\Lbac
B\left(f^{m}(s),\left|\frac{(1-\delta')l_m^{n-1}\pi}{2Q'(f^m(s))}\right|\right)$ in the sense that as $m\to\infty$ the measure of density of the complement of $f^{m-m_s}(A_{m,n})$ in this annulus and vice versa tends to zero.
This annulus is contained in $G_\infty$ and
its diameter is more than twice as large as the width of the gaps in between the channels of $\Gamma$. Thus it has to intersect these channels. More precise, the diameter assures that the density of $\Gamma$ in $f^{m-m_s}(A_{m,n})$ is bounded below by some constant $c_2$. Since the distortion of $f^{m-m_s}$ on $A_{m,n}$ is arbitrarily close to one for large $m$, for $c_3<c_2$ this carries over by (\ref{qs2}) to \begin{eqnarray}\frac{\meas(A_{m,n}\cap \Gamma)}{\meas(A_{m,n})}&>& c_3. \label{c3} \end{eqnarray} \\
\begin{figure}
\begin{center}
\scalebox{0.5}{\input{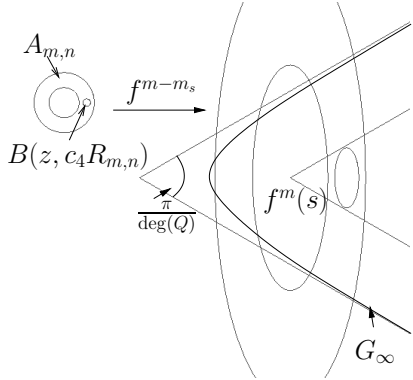}}
\end{center}
\caption{Constructions to obtain estimates for the density of $\Gamma$ in $A_{m,n}$}\label{smalllargen}
\end{figure}
For larger $n$ the distortion is still bounded by $K:=K_{1/\sqrt{2}}$, such that one could call $f^{m-m_s}(A_{m,n})$ a
$K$-quasi-annulus, whose center $f^{m}(s)$ lies in $G_{\infty}$ and whose diameter is far larger than the gaps in between
the channels of $\Gamma$. For $0<c_4<\sin\left(\frac{\pi}{2\deg(Q)}\right)/K$ the following is easy to verify: For large $m$ there exists $z\in \partial B(f^{m_s}(s),\frac{R_{m,n}+R_{m,n-1}}{2})$, such that $f^{m-m_s}(B(z, c_4 R_{m,n}))  \subset G_\infty$. Here the choice of $c_4$ guarantees that $f^{m-m_s}(B(z, c_4 R_{m,n}))$ lies in a sector of angle $\pi/\deg(Q)$ and corner $f^m(s)$. The boundary of the component of $G_\infty$ containing $f^m(s)$ is tangent to the boundary of a sector of the same angle and corner $0$. Thus it is sufficient to choose $z$, maximizing the distance of $f^{m-m_s}(z)$ to $\C\Lbac\G_{\infty}$. This is displayed in figure \ref{smalllargen}. We note that $c_4<\frac{1}{4}$, such that $B(z, c_4 R_{m,n})\subset
A_{m,n}$ for large $m$. Using the family $\F_\infty$ from above as a cover, it is easy to see that
the density of $\Gamma$ in $f^{m-m_s}(B(z, c_4 R_{m,n} ) )$ is bounded below by any
$0<c_5<c'$ and $m$ large enough. This carries over by (\ref{qs2}) to $B(z,
c_4\diam(f^{m-m_s}(A_{m,n})))$, in which the density of $\Gamma$ is at least $\frac{c_5}{K^2}$. We assume  $c_3\le\frac{c_5 c_4^2}{K^2}$ and get (\ref{c3}) for all $n$ and $m$ if $m_0$ is large enough. This
carries over by (\ref{qs2}) to the elements of $\F_s$ and completes the proof for $c=\frac{c_3}{C_2^2C_3^2}$.\qed
\begin{lem} \label{AsubsetB}Let $B\subset\C$ be finite, such that every $b\in B$ escapes exponentially. Suppose that every singularity of the inverse is either pre-periodic, escapes exponentially or is contained in some attractive basin. Suppose further that $A\not\subset O^+(B)$. Then the set $\{z\in J(f):\omega(z)\subset O^+(B)\cup\{\infty\}\}$ has zero measure.
\end{lem}
\proof We assume positive measure of the set above. Then for some $\ep>0$ also the set $X:=\{z\in J(f):\{\infty\}\not=\omega(z)\subset O^+(B)\cup\{\infty\},\forall a\in A\Lbac O^+(B): \dist(O^+(z),a)>\ep \}$ has positive measure. We may assume $\omega\{z\}\not=\{\infty\}$ since A. Eremenko and M. Lyubich \cite{ere} proved that the
set of escaping points $I(f)$ has measure zero. The assumption $\dist(O^+(z),a)>\ep$ is permissible since $X$ is contained in the countable union of the sets $\{z :\forall a\in A :\dist(a,O^+(z))>\ep_n)\}>0$ with $\ep_n\to 0$.\\ Since every $b\in B$ escapes exponentially, there exists some $\delta>0$
with $|f^{n+1}(b)|\ge \exp(|f^n(b)|^{4 \delta})$ for every $b\in B$ and $n$ large enough. Let $z_0$ be
a density point of $X$. Since the iterates of $z_0$ do not tend to infinity, there exists a convergent subsequence $f^{\beta(n)}(z_0)$, whose limit must be of the form $f^{n_0}(b)$ with $n_0\in \N\cup\{0\}$ and $b\in B$. We may assume
$n_0\ge m_b:=\max(\{m\in\N:f'(f^{m-1}(b))=0\}\cup\{0\}) $. For all $n\in\N$ we define $\alpha(n)\ge \beta(n)$ smallest possible with
\begin{eqnarray}
|f^{\alpha(n)}(z_0)-f^{\alpha(n)-\beta(n)+n_0}(b)|&\ge&
|f^{\alpha(n)}(z_0)|^{1-\deg(Q)+3\delta} \end{eqnarray}
and $B_n:= B(f^{\alpha(n)}(z_0), |f^{\alpha(n)}(z_0)|^{1-\deg(Q)+2\delta})$.
We will see that for large $n$ the inverse branch $g_n$ of $f^{\alpha(n)}$, mapping $f^{\alpha(n)}(z_0)$ to $z_0$, may be extended with uniformly bounded distortion to $B_n$. Furthermore we show that the density of $X$ in $B_n$ does not tend to one. This carries over to $g_n(B_n)$ by (\ref{qs2}). Finally we show $\diam(g_n(B_n))\to 0$, being a contradiction to the choice of $z_0$ as a density point of $X$.\\
Since $\sing(f^{-1})$ is bounded, we can extend every branch of $f^{-1}$ to $B(z,(1-\delta')|z|)$ for every $\delta'>0$, and $|z|$ large enough. Thus we can extend the branch of $f^{-1}$, mapping $f^{\alpha(n)}(z_0)$ to  $f^{\alpha(n)-1}(z_0)$, to $B_n$, such that the distortion is bounded by a constant, which can be chosen arbitrarily close to one for $n$ sufficiently large, and its image is contained in $B(f^{\alpha(n)-1}(z_0), |f^{\alpha(n)}(z_0)|^{-\deg(Q)+3\delta})$. Due to the definition of $\alpha$ the last set is itself contained in $B(f^{n_0+\alpha(n)-\beta(n)-1}(b),2|f^{n_0+\alpha(n)-\beta(n)-1}(b)|^{1-\deg(Q)+3\delta} )$. The choice of $n_0\ge m_b$ and the same argument as above applied $\alpha(n)-\beta(n)-1$ times allows us to extend the branch of $(f^{\alpha(n)-\beta(n)-1})^{-1}$, mapping $f^{n_0+\alpha(n)-\beta(n)-1}(b)$ to  $f^{n_0}(b)$, to the set $B(f^{n_0+\alpha(n)-\beta(n)-1}(b),(1-\delta')|f^{n_0+\alpha(n)-\beta(n)-1}(b)|) $. Again this implies that its distortion on $B(f^{n_0+\alpha(n)-\beta(n)-1}(b),2|f^{n_0+\alpha(n)-\beta(n)-1}(b)|^{1-\deg(Q)+3\delta} )$ tends to $1$ as $n$ tends to infinity.  Moreover the pre-image of $B_n$ under $f$, which is contained in $B(f^{\alpha(n)-1}(z_0), |f^{\alpha(n)(z_0)}|^{-\deg(Q)+3\delta})$ is mapped to $B(f^{\beta(n)}(z_0),r_n)$ with $r_n/|f^{\beta(n)}(z_0)-f^{n_0}(b) |\to 0$. Since $P(f)$ does not accumulate at $f^{n_0}(b)$, we can extend the branch of the inverse of $f^{\beta(n)}$, mapping $f^{\beta(n)}(z_0) $ to $z_0$, to $ B(f^{\beta(n)}(z_0),|f^{\beta(n)}(z_0)-f^{n_0}(b) | )$. Thus $g_n$ exists and its distortion tends to one as $n$ tends to infinity. \\
For $s\in A\cup\{\infty\}$ we define $G_s$ as in lemma \ref{dengamma}. Due to our choice of $\diam(B_n)$, the density of
$\C\Lbac G$ in $B_n$ tends to zero as  $n$ tends to infinity.
For $M$ large enough and $s\in A\Lbac O^+(B)$ we
have $G_s\cap X=\emptyset$.  For $s\in A\cap O^+(B)\cup\{\infty\}$ we can apply
lemma \ref{dengamma} to $\Gamma:=\C\Lbac X$ and obtain a family $\F$ of disjoint domains covering $B_n\cap
G_s$ up to measure zero, such that the density of $\Gamma $ in all of these is bounded below by some positive constant $c$. The diameter of these domains is much smaller than the diameter of $B_n$, such that we can neglect the ones intersecting the boundary of  $B_n$. Thus the density of $X$ in $B_n$ does not tend to one. \\
It remains to show that $\diam(g_n(B_n))\to 0$, which is equivalent to $\left|(f^{\alpha(n)})'(z_0)\right||f^{\alpha(n)}(z_0)|^{\deg(Q)-1-2\delta}\to \infty$.\\ First we show that for $\ep'>0$ small enough, $b\in A\cap O^+(B)$, $z\in B(s,\ep')$ and \begin{eqnarray}\label{mz}&m:=& \min\left\{ m\in\N:|f^{m}(z)-f^{m}(s)|\ge |f^{m}(s)|^{1-\deg(Q)+3\delta}\right\}\end{eqnarray} we have \begin{eqnarray}\label{neunviertela}|(f^m)'(z)|&=& \frac{2|f^m(z)|^{1-\deg(Q)+\frac{9\delta}{4}}}{|z-s|}.\end{eqnarray}
To see this, we assume that $\ep'$ is small enough to guarantee that $f^{m-1}$ maps  $B(s,|z-s|)$ to $B(f^{m-1}(s),2|f^{m-1}(s)|^{1-\deg(Q)+3\delta} )$ with a distortion that is very close to one. This may be achieved since for $\ep'$ small $m$ is large and, with lemma \ref{finj}, one can show that the function $f^{m-1}$ is injective on a region far larger than $B(s,|z-s|)$. Thus (\ref{estf}) implies that the magnitude of $\mbox{Re}(Q(f^{m-1}(x)))$ is of the same for all  $x\in B(s,|z-s|)$. With (\ref{estf}) this carries over to $|f^m(x)|$. In particular we may assume that \begin{eqnarray}\label{viertel}|f^{m}(s)|^{1-\frac{\delta}{4}}\le&|f^m(x)|&\le \;|f^{m}(s)|^{1+\frac{\delta}{4}}.\end{eqnarray} Then (\ref{neunviertela}) emerges as follows 
\begin{eqnarray}|(f^m)'(z)|&=&| f'(f^{m-1}(z))||(f^{m-1})'(z)|\nonumber\\
&\overset{}{\ge}& |  f^m(z)|| Q'(f^{m-1}(z))||(f^{m-1})'(z)|\nonumber\\&\overset{}{\ge}&  \sup_{x\in [s,z]}|f^m (x)||f^{m}(s)|^{-\frac{\delta}{2}}| Q'(f^{m-1}(z))||(f^{m-1})'(z)| \nonumber\\
&\overset{}{\ge}& \frac{1}{2} \sup_{x\in [s,z]}|f^m (x)||f^{m}(s)|^{-\frac{\delta}{2}}| Q'(f^{m-1}(x))||(f^{m-1})'(z)| \nonumber\\
&\overset{}{\ge}&\frac{1}{2} \sup_{x\in [s,z]}|(f^m)' (x)| |f^{m}(s)|^{-\frac{\delta}{2}}| \nonumber\\
&\overset{}{\ge}& \frac{|f^m(z)-f^{m}(s)| |f^{m}(s)|^{-\frac{\delta}{2}}}{2|z-s|}\nonumber\\
&\overset{}{\ge}& \frac{|f^{m}(s)|^{1-\deg(Q)+\frac{5\delta}{2}}}{2|z-s|}\nonumber\\
&\overset{}{\ge}&  \frac{2|f^m(z)|^{1-\deg(Q)+\frac{9\delta}{4}}}{|z-s|}.\nonumber\end{eqnarray}
Here the second and fifth estimate follows with (a), while the third and last follow with (\ref{viertel}) and the mean value theorem  and the Definition of $m$ imply the sixth and seventh estimate respectively.
From (\ref{estf}) follows the existence of $R>0$, such that $\frac{1}{2}|f'(z)| \le\frac{|f(z)-\overline{s}(z)|}{|Q'(z)|}\le 2|f'(z)| $ for any $z\in\C$ with $|z|\ge R$. Now let $I$ be the set of $j\in\N$, such that $f^j(z_0)\not\in B(0,R)$ and $|f'(f^j(z_0))|\le \frac{\ep'}{2}|Q'((f^j(s))|$. Then for $j\in I$ there exists $s\in A$ with
\begin{eqnarray}\label{j+1}|f^{j+1}(z_0)-s|&\le& \frac{2|f'(f^j(s))|}{|Q'((f^j(s))|},\end{eqnarray} which is at most $\ep'$. We assume $\ep'<\ep$, such that $s\in O^+(B)$. We choose $m$ as in (\ref{mz}) for $z=f^{j+1}(s)$. Then we have $m\le\alpha(n)$ and (\ref{neunviertela}) together with (\ref{j+1}) imply that
\begin{eqnarray}|(f^{m})'(f^j(z_0))|=|f'(f^j(z_0)) (f^{m-1})'(f^{j+1}(z_0) )|\ge |Q'(f^j(z_0))|
|f^{m+j}(z_0)|^{1-\deg(Q)+\frac{9\delta}{4}}.\label{neunviertel}
\end{eqnarray}
Since $\omega(z_0)\subset \overline{O^+(B)}$, there exists $j_0\in\N$, such that for all $j\ge j_0$ the point $f^j(z_0)$ is not contained in the compact set $\overline{B(0,R)}\Lbac\left( \bigcup_{n\in\N,b\in B}f^n(B(b,\frac{\ep'}{4})\right)$. In the case $\deg(Q) > 1$ we assume $|Q'(z)|\ge \frac{2}{\ep'}$
and in the case $\deg(Q)=1$ we assume $R\ge \frac{2}{|q|}+\max_{a\in A}|a|$, where $q$ is the leading coefficient of $Q$.
Then for all $j_0\le j\not\in I$ with $|f'(f^j(z_0))|\le 1$ we have $f^j(z_0)\in f^n(B(b,\frac{\ep'}{4}))$ for some $n$ and $b\in B$. The exponential escape of $b$ implies the existence of some $m\in\N$ only depending on $n$ and $b$, for which $|(f^{m})'(f^j(z_0))|\ge 1$. \\
Let $\Delta_j:=\min\{i\in I\cup\{\alpha(n)\}:i>j \}-m_j-j$. We assume $j_0\in I$ and $j\not\in I$ for $j<j_0$. (Otherwise we define $j_0\le j_0'\in \N\cap I$ smallest possible and $I':=I\Lbac \{1,..j_0'\}$ and continue with those.) With the chain rule we get \begin{eqnarray}\left|\frac{(f^{\alpha(n)})'(z_0)}{f^{\alpha(n)}(z_0)|^{1-\deg(Q)+2\delta}}\right|&\!\!\!\!\!=&\!\!\!\!\!\left|\frac{(f^{j_0})'(z_0)}{f^{\alpha(n)}(z_0)^{1-\deg(Q)+2\delta}}\right|\!\!\!\!\!\!\!\!\!\!\!\!\!\!\!\!\prod_{\;\;\;\;\;\;\;\;\;\;\;\;\;\;\;j\in I,j_0\le j<\alpha(n)} \!\!\!\!\!\!\!\!\!\!\!\!\!\!\!\!|(f^{\Delta_j })'(f^{m_j+j}(z_0))(f^{m_j})'(f^j(z_0)) | \nonumber\end{eqnarray}
With (\ref{neunviertel}) we obtain a lower estimate of the product above, in which most of the factors in the product cancel out each other. More precisely for each $j\in I$ except the first and the last in the product the factor $ |f^{m_j+j}(z_0)|$ only remains with a power of $\frac{9}{4}$. If $\Delta_j=0$ this is follows directly from (\ref{neunviertel}) by considering the $j$-th and the $j+1$-st factor of the product together. Of course there appears the factor $|q|\deg(Q)$ as the leading coefficient of $|Q'|$.  
If $\Delta_j\not=0$, we have $m_j+j\not\in I$, such that $|(f^{\Delta_j })'(f^{m_j+j}(z_0))|\ge |f'(f^{m_j+j}(z_0))|\ge \frac{\ep'}{2}|Q'(f^{m_j+j}(z_0))|$. This implies the same as above with the factor $\frac{\ep'|q|^2}{2}$, which we assume to be smaller than $|q|\deg(Q)$.\\
Finally we know due to the definition of $m_j$, that for the last $j$ in the product we have $j+m_j=\alpha(n)$ such that this factor cancels out with the denominator in front of the product up to $|f^{\alpha(n)}(z_0)^{\frac{\delta}{4}}|$ and we get
\begin{eqnarray}\left|\frac{(f^{\alpha(n)})'(z_0)}{f^{\alpha(n)}(z_0)|^{1-\deg(Q)+2\delta}}\right|&\!\!\!\!\!\overset{}{\ge}& \!\!\!\!\!|f^{\alpha(n)}(z_0)|^{\frac{\delta}{4}}|(f^{j_0})'(z_0)||Q'(f^{j_0}(z_0))|\!\!\!\!\!\!\!\prod_{j\in I,j_0<j<\alpha(n)-m_j} 
\!\!\!\!\!\!\!\!\!\!\!\!\!\!\frac{\ep'q^2}{2} |f^{m_j+j}(z_0)|^{\frac{9\delta}{4}}, \nonumber\nonumber
\end{eqnarray} 
which tends to infinity as $n$ does so, and thus completes the proof. \qed
With this lemma we can finally prove the theorems \ref{andereinclusion} and \ref{class}.
\proofof{theorem \ref{andereinclusion}}
The assumptions on
the singular orbits guarantee $P(f)'\cap J(f)=\emptyset$. An indifferent periodic point in the Julia set must be an accumulation point of $P(f)$. This is a well known fact, whose proof may be found in \cite{milnor}, where it is stated for rational functions. However only minor changes are necessary for the transcendental case. Therefore all periodic points in $J(f)$ are repelling. Due
to theorem \ref{asatz}, we have $\omega(z)\subset P(f)$ for almost every $z\in J(f)$. If $O^+(z)$ accumulates at a
repelling periodic point, $\omega(z)$ accumulates also at this point. This follows from the fact that $O^+(z)$ accumulates at every compact annulus $\{z: r \le |z-p|\le 2|(f^n)'(p)|r \}$, if $r>0$ small enough and $n$ is the period of the repelling periodic point $p$. Thus for almost every $z\in J(f)$ we have
$\omega(z)\subset \overline{O^{+}(B)}$, if $B$ is the set of singularities, that escape exponentially. Now lemma \ref{AsubsetB} implies that the set of points, that do not accumulate at all asymptotic values, has measure zero. This concludes the proof for the inclusion $\omega(z)\supset A$. \\
Now assume that there exists some point $s\in B\Lbac \overline{O^+(A)}$, such that the set $X':=\{z\in J(f):s\in \omega(z)\subset \overline{O^+(B)}\}$ has positive measure. Then the whole proof of lemma \ref{AsubsetB} works identically, with $X'$ instead of $X$ and with the only difference, that at the point, where lemma \ref{dengamma} is used, we now use the measure estimate of theorem \ref{asatz} instead. More precise, instead of using the family $\F$ from lemma \ref{dengamma}, to see that the density of $\Gamma$ in $B_n:=B(f^{\alpha(n)(z_0},|f^{\alpha(n)(z_0}|^{1-\deg(Q)+2\delta})$ is bounded below, we argue as follows, to see that the density of $T(f):=\{z:\omega(z)\subset \overline{O^+(A)}\}$ in $B_n$ is bounded below. Since $O^+(T(f))$ is disjoint from $O^+(X')$, we get again that $X'$ contains no density point, contradicting the assumption of positive measure.\\
To obtain the new family, we proceed as in the proof for $\meas(F(f))<\infty$ of theorem \ref{gsatz}. We cover the set $\{z:\dist(z,\C \Lbac G)\ge
2|z|^{1-\deg(Q)+\delta}\}$ with a family $\S$ of squares $S\subset \{z:\dist(z,\C \Lbac G)\ge |z|^{1-\deg(Q)+\delta}\}$ with
$\sup_{z\in S}|z|^{1-\deg(Q)}\le\diam(S) \le4\sup_{z\in S}|z|^{1-\deg(Q)}$. The density of $T(f)$ in any $B_n\supset S\in\S$ is, due to theorem \ref{asatz}, very close to one for large $n$. In particular it is bounded below by some positive constant $c$. The diameter of $B_n$ implies that the density of $\{z:\dist(z,\C \Lbac G)\le
2|z|^{1-\deg(Q)+\delta}\}$ tends to zero with $n$. The same is true for the union of those squares in $\S$, that intersect the boundary of $B_n$. This gives the estimate needed to proceed with the proof of lemma \ref{AsubsetB}.\qed 
\proofof{theorem \ref{class}} If all asymptotic values escape exponentially we can apply theorem (\ref{gsatz}) and obtain a set of positive measure, whose orbits accumulate only at the orbits of the asymptotic values and the point infinity. In particular the function is not recurrent. We assume now that the set of pre-periodic asymptotic values is non-empty and the post-critical
case holds. Since $P(f)'\cap J(f)=\emptyset$, there are again no indifferent periodic points. From lemma \ref{AsubsetB} we know that the
orbit of almost every $z\in J(f)$ accumulates at least at one point in $P(f)$, which does not escape exponentially and
 thus has to be pre-periodic. By continuity $O^+(z)$ accumulates at a repelling periodic point. As above this implies that $\omega(z)$ accumulates at this repelling periodic point.
 Since $P(f)$ has no cluster points in $\C$, this is a contradiction. Thus (i) of theorem \ref{bockclass} is satisfied and $f$ is recurrent and ergodic on $J(f)=\C$. \qed

\section{Other applications}
As mentioned in the remark \ref{remk}, one can use the theorem \ref{asatz}, in order
to obtain positive measure for the escaping set $I(f)$. As one example one can consider the following family
containing the sine and cosine family, for which this result was proved by C. McMullen in \cite{mcmullen}.
\begin{thm}\label{if}
Let $f(z):=P(z)\exp(Q(z))+ \tilde{P}(z)\exp(\tilde{Q}(z))$ for polynomials $P,\tilde{P}\not=0$ and $Q,\tilde{Q}$, such
that $n:=\deg(\tilde{Q})=\deg(Q)\ge 0$ and the arguments of their $n-th$ coefficients $q,\tilde{q}$ differ by some odd multiple of
$\frac{\pi}{n}$. Then $ \meas(I(f))>0$. If $\tilde{Q}=-Q$ and $n\ge3$ then $\meas(\C\Lbac I(f))<\infty$.\end{thm}
\sketch As in \ref{gsatz} on can show, that for any
  $0<\delta<\beta< 1$,
 $-1<\delta_1<n-1<\delta_2$, $M$ large enough, $A:=\emptyset$ and $G:=\{z:|\arg(z)-\frac{(2k+1)\pi-2\arg(q)}{2 n}|\le |z|^{\delta-1}\}$
the conditions (a) and (b) of theorem \ref{asatz} are satisfied, while condition (c) is trivial. The theorem implies the
first part. The second follows as in the proof of \ref{gsatz} choosing $1-n<\delta<\beta<-1$.
 \qed
As an example we consider the function \[ f(z)\;:=\;\exp(z^3)-\exp(-z^3).\]
Its Fatou set is not empty, since it contains a super attractive basin around zero. The theorem above gives however $0<\meas(\C\Lbac I(f))<\infty$. In the figure below the Fatou set is black. The picture shows the part of the plane given by $\{z:|\mbox{Re}(z)|\le 2,|\mbox{Im}(z)|\le 2 \}$.\\
\begin{figure}[!ht!]
\begin{center}
\includegraphics[width=50mm]{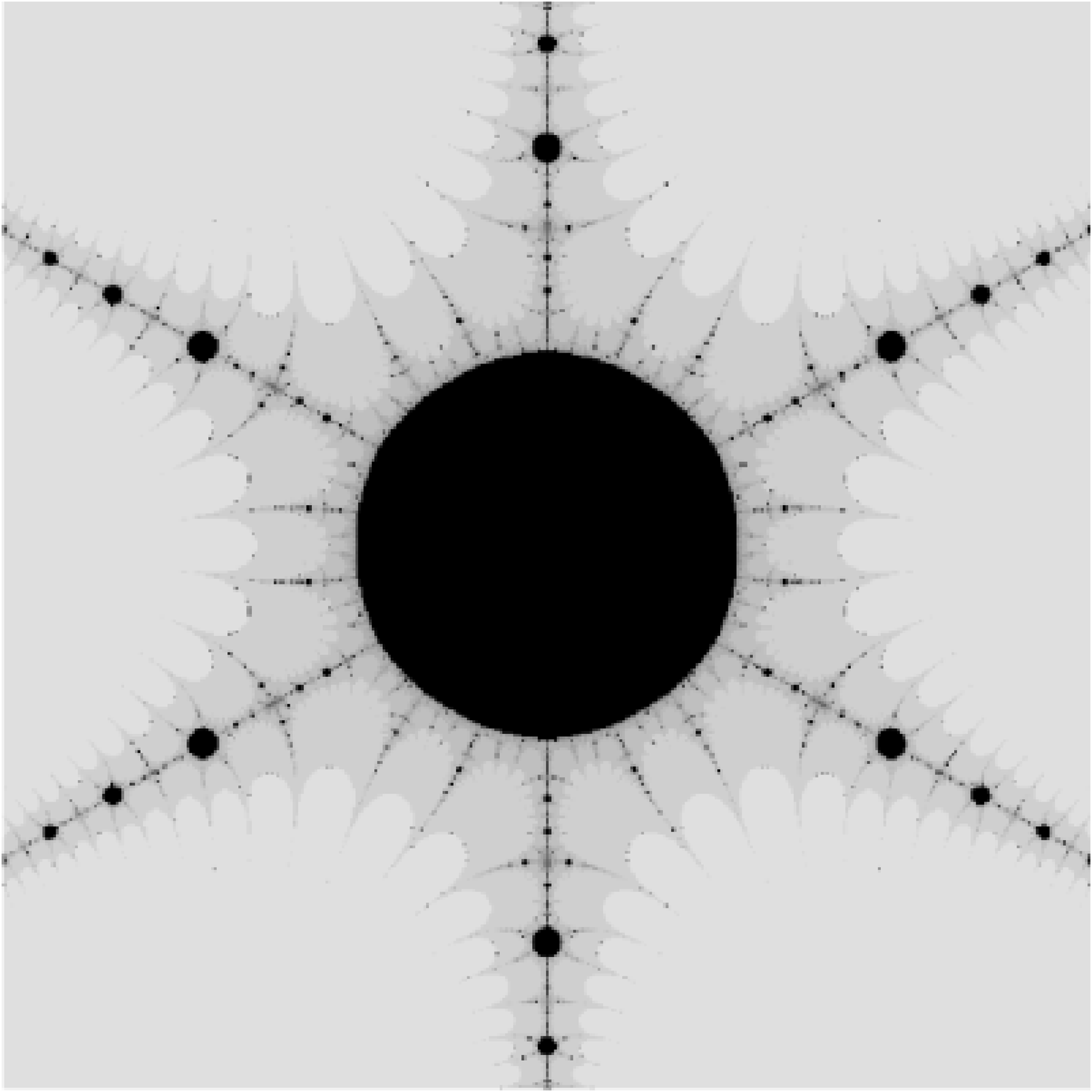}
\caption{The Fatou set of $f(z)=\exp(z^3)-\exp(-z^3)$}
\end{center}
\end{figure}
 The functions discussed in the last chapter have rational Schwarzian
 derivative $S(f):=\frac{f'''}{f'}-\frac{3}{2}\left(\frac{f''}{f'} \right)^2$. 
There are many things known about functions, whose Schwarzian derivative is a polynomial. The asymptotic behavior of functions with this property has already been studied by E. Hille \cite{hille} and R. Nevanlinna \cite{nevanlinna1}. Most things
 carry over to the rational case, which has been studied by G. Elfving \cite{elfving}.  
It is easy to see that a critical point of $f$ is a pole of $S(f)$. Thus functions with a rational Schwarzian derivative have only finitely many critical points. If $S(f)(z)=c z^n(1+\oo(1))$ as $z\to\infty$ with $c\not=0$ and $n\ge 0$, there are $n+2 $ so called critical  rays  defined by $\arg z=\phi$ with $\arg c +(n+2)\phi =0 (\mod 2\pi)$. It turns out that these divide the complex plane into $n+2$ sectors in which the asymptotic behavior of $f$ is known very well. If $z$ tends to infinity in a non-critical direction, $f$ tends to an asymptotic value, which is the same for all directions inside the same sector. Thus $f$ has only finitely many asymptotic values. Similar as in the proof of theorem \ref{gsatz} one can show, that the conditions (a) and (b) of theorem \ref{asatz} are satisfied. If one of these asymptotic values happens to be $\infty$, points and also asymptotic values may escape exponentially inside the corresponding sector satisfying condition (c). However these functions may have
 infinitely many poles, such that points can also escape exponentially, ``jumping from pole to pole'',
 without satisfying the condition (c) of theorem \ref{asatz}. The poles are however contained in small neighborhoods around these critical rays. Thus we can formulate another more geometric condition in order to guarantee
 condition (c). More precisely we get
\begin{thm}
Let $f$ be a meromorphic function with rational Schwarzian derivative, whose behavior at infinity is of the form $c z^n(1+\oo(1))$ with $c\not=0$ and $n\ge-1$. Suppose that all asymptotic
values $s$ tend to $\infty$ under iteration and there exists some $\ep>0$, such that $|\arg(f^m(s))-\frac{2\pi k +\arg(c)}{n+2}|\ge |f^m(s)|^{\ep-\frac{n+2}{2}}$ for almost all $m\in\N$ and all $k\in\{0,1,...,n+1\}$. Then $\meas(J(f))>0$ and $\omega(z)\subset P(f)$ for almost every $z\in J(f)$. If $n\ge 3 $ it follows that $\meas(F(f))<\infty$.\end{thm}
\sketch
The principle is exactly as the proof of theorem \ref{gsatz}. First one has to check that the properties of theorem \ref{asatz} are satisfied. This gives us measure estimates of $T(f)$ that imply case (ii) of theorem \ref{bockclass}. We obtain $T(f)\subset J(f)$ again from the absence of Baker and wandering domains, which once more follows from the finiteness of $\sing (f^{-1})$ (see remark \ref{remk}). For meromorphic functions with polynomial Schwarzian derivative this has also been shown by R. L. Devaney and L. Keen in \cite{devaneykeen}.\\
To check the properties we briefly summarize how to obtain estimates of the asymptotic behavior of $f$. We refer to the post-graduate notes of Jim Langley \cite{langley} for more details.
It is easy to see that functions with the rational Schwarzian derivative $S(f)=2 A$ coincide with those quotients $f_1/f_2$ of two linearly independent solutions of the differential equation
$f_i''+Af_i=0$. Moreover the asymptotic behavior of these solutions is known very well. For a critical ray with argument $\phi$ and $R_0>0$ large we define 
\begin{eqnarray}Z(z):=\int_{2 R_0e^{i\phi}}^z A(t)^{1/2}=\frac{2c^{1/2}}{n+2}z^{(n+2)/2}\left(1+\OO\left(\frac{\ln|z|}{|z|}\right)\right),\;\;\text{for}\; z\to\infty\end{eqnarray}
in the set $\{z:R_0\le |z|,|\arg z -\phi|\le \frac{2\pi}{n+2}\}$. Then it is easy to see that for $\delta'>0$ and $R_1$ large enough $Z$ is univalent in the set $S_1:=\{z:|z|\ge R_1,\; |\arg(z)-\phi|<\frac{2\pi}{n+2}-\delta'\}$.
With the Liouville transformation $ W_i(Z)=A(z)^{1/4} f_i(z)$ and $F_0(Z):=A''(z)/4A(z)^2-5A'(z)^2/16A(z)^3$ we get 
\begin{eqnarray}
\frac{ \partial^2 W_i }{ \partial Z^2 }+(1-F_0(Z))W_i&=&0.
\end{eqnarray}
This equation has been integrated asymptotically by Hille \cite{hille} and his method has been used by many others afterwards. The following theorem may be found explicitly in \cite{langley}.
\begin{thm}[Hille, Langley]
Let $c'>0$ and $0<\ep'<\pi$. Then there exists a constant $d'>0$, depending only on $c$ and $\ep'$, with the following properties.
Suppose that $F$ is analytic, with $|F(z)|\le c'|z|^{-2}$, in
\begin{eqnarray}\Omega:=\{z:1\le R_0 \le |z|\le R_1<\infty,|\arg z |\le\pi-\ep'\}.\nonumber
\end{eqnarray}
Then the equation 
\begin{eqnarray}
\omega''+(1-F(z))\omega &=&0
\end{eqnarray}
has two linearly independent solutions $U,V$ satisfying
\begin{eqnarray}
U(z)=\exp(-iz)(1+\delta_1(z))&,& U'(z)=-i\exp(-iz)(1+\delta_2(z)),\nonumber\\
V(z)=\exp(iz)(1+\delta_3(z))&,& V'(z)=i\exp(iz)(1+\delta_4(z)),
\end{eqnarray}
such that $|\delta_i(z)|\le d'|z|^{-1}$ for  $z\in \Omega\Lbac\{z:\mbox{Re}(z)<0,|\mbox{Im}(z)|<R\}.$
\end{thm}
\remark $\Omega$ may be replaced by \begin{eqnarray}\Omega':=\{z:1\le R_0\le |z|\le R_1<\infty,|\arg z -\pi|\le\pi-\ep'\}\nonumber
\end{eqnarray} and also by the unbounded region \begin{eqnarray}\Omega'':=\{z:1\le R_0\le |z|< \infty,|\arg z |\le\pi-\ep'\}.\nonumber
\end{eqnarray}
To see this, we take a sequence $R_k\to\infty$ and obtain solutions $U_k$,$V_k$ with uniformly bounded $\delta_{i,k}$ in $\Omega_k$, where $\Omega_k$ is $\Omega$ with $R_1$ replaced by $R_k$. Therefore both form a normal family, and a subsequence of $U_k$,$V_k$ converges in $\Omega''=\bigcup_{k\in\N}\Omega_k$.\\\quad\\
Thus for every $j\in \{1,..,n+2\}$ and every critical ray with argument $\phi_j$ there are constants $a_j,b_j,c_j,d_j\in\C$, such that \begin{eqnarray}
f(z)&=& \frac{a_j U(Z )+b_j V(Z)}{c_j U(Z )+d_j V(Z)}\\&\!\!\!\!\!\!\!\!\!\!\!\!\!\!\!\!\!\!\!\!\!\!\!\!\!\!\!\!\!\!\!\!\!\!\!\!=&\!\!\!\!\!\!\!\!\!\!\!\!\!\!\!\!\!\!\!\!\!\frac{\left(a_j\exp\left(\frac{-2ic^{1/2}}{n+2}z^{\frac{n+2}{2}}(1+\OO(\frac{\ln|z|}{|z|}))\right) +b_j\exp\left(\frac{2ic^{1/2}}{n+2}z^{\frac{n+2}{2}}(1+\OO(\frac{\ln|z|}{|z|}))\right)\right)(1+\OO(|z|^{\frac{-1}{2}}))   }{\left( c_j\exp\left(\frac{-2ic^{1/2}}{n+2}z^{\frac{n+2}{2}}(1+\OO(\frac{\ln|z|}{|z|}))\right)+d_j\exp\left(\frac{2ic^{1/2}}{n+2}z^{\frac{n+2}{2}}(1+\OO(\frac{\ln|z|}{|z|}))\right)\right)(1+\OO(|z|^{\frac{-1}{2}}))}\nonumber
\end{eqnarray}
for $z\to\infty$ in $S_j:=\{z:1\le R_0\le |z|,\; |\arg(z)-\phi_j|<\frac{2\pi}{n+2}-\delta\}$.
Thus $f$ tends to $a_j/c_j$ in $S_j^+:=\{z\in S_j:\arg(z)>\phi_j|\}$ and $f$ tends to $b_j/d_j$ in $S_j^-:=\{z\in S_j: \arg(z)<\phi_j\}$. If $c_j$ or $d_j$ happen to be zero, while $a_j$ or $b_j$ are not, we obtain a sector, on which $f$ tends to $\infty$, such that points may escape exponentially in this sector.
We get a similar estimate for the derivative, such that we can prove with similar arguments as in the proof of theorem \ref{gsatz} that $f$ satisfies the conditions of theorem \ref{asatz} for the choices $0<\delta< \ep$, $\delta-\frac{n}{2}<\beta<1$, $-1<\delta_1<\frac{n}{2}<\delta_2$, $M$ large enough, and $G:=\bigcup_{1\le j\le n+2} \{z\in S_j:|\mbox{Im}(Z_j(z))|\ge |Z_j(z)|^{\frac{2\delta}{n+2}}\} $, where $Z_j$ is the upper change of coordinates $Z$ for the sector $S_j$. If $n\ge 3$, we can choose $\delta<\frac{n-2}{2}$. Then the proof for $\meas(F(f))<\infty$  works just as in theorem \ref{gsatz}. \qed

\bibliographystyle{acm}
\bibliography{dipl}
\end{document}